\documentclass[letterpaper,12pt]{article}
\usepackage{amssymb}
\usepackage{amsmath}
\usepackage{makeidx}
\usepackage{amscd}
\usepackage[arrow,matrix]{xy}

\newtheorem{theo}{Theorem}[section]
\newtheorem{prop}[theo]{Proposition}
\newtheorem{lem}[theo]{Lemma}
\newtheorem{cor}[theo]{Corollary}
\newtheorem{defi}[theo]{Definition}

\def \deg {{\rm{deg}}}
\def \Br {{\rm{Br}}}

\def \Ga {{\Gamma}}

\def \Pic {{\rm {Pic}}}

\def \Gal {{\rm{Gal}}}
\def \Ker {{\rm{Ker}}}
\def \Im {{\rm {Im}}}

\def \A{{\mathbb A}}
\def \P{{\mathbb P}}
\def \Spec {{\rm{Spec}}}
\def \dim {{\rm{dim}}}
\def \Hom {{\rm {Hom}}}
\def \End {{\rm {End}}}
\def \Pic {{\rm {Pic}}}
\def \GL {{\rm {GL}}}

\def \Aut{{\rm Aut}}

\def\ov{\overline}

\def \Z {{\mathbb Z}}
\def \Q {{\mathbb Q}}
\def \F {{\mathbb F}}

\def \Id {{\rm Id}}

\def \rk {{\rm{rk}}}

\def\G{{\mathbb G}}

\def\GG{{\cal G}}
\def\HH{{\cal H}}
\def\C{{\mathbb C}}

\def\T{{\cal T}}

\def\lra{\longrightarrow}

\def\H{{\rm H}}

\def\sA{{\cal A}}

\def\L{{\cal L }}

\def\NS{{\rm NS\,}}

\def\O{{\cal O}}
\def\L{{\cal L}}

\def\codim{{\rm codim}}

\def\char{{\rm char}}
\def\Ga{\Gamma}

\def\et{{\rm{\acute et}}}

\def\bS{{\bf S}}
\def\bA{{\bf A}}

\newcommand{\bthe}{\begin{theo}}
\newcommand{\ble}{\begin{lem}}
\newcommand{\bpr}{\begin{prop}}
\newcommand{\bco}{\begin{cor}}
\newcommand{\bde}{\begin{defi}}
\newcommand{\ethe}{\end{theo}}
\newcommand{\ele}{\end{lem}}
\newcommand{\epr}{\end{prop}}
\newcommand{\eco}{\end{cor}}
\newcommand{\ede}{\end{defi}}

\title{Kummer varieties and their Brauer groups}
\author{Alexei N. Skorobogatov and Yuri G. Zarhin}
\date{\today}
\begin{document}
\baselineskip=15pt
\maketitle

\centerline{\em{to Yuri Ivanovich Manin on his 80th birthday, with admiration and gratitude}}

\begin{abstract}
\noindent We study Kummer varieties attached to 2-coverings of abelian varieties
of arbitrary dimension. Over a number field we show that the subgroup 
of odd order elements of the Brauer group
does not obstruct the Hasse principle. Sufficient conditions for the triviality
of the Brauer group are given, which allow us to give an example of
a Kummer K3 surface of geometric Picard rank 17 over the rationals
with trivial Brauer group. We establish the non-emptyness of the Brauer--Manin set 
of everywhere locally soluble Kummer varieties attached
to 2-coverings of products of hyperelliptic Jacobians with large Galois action on 2-torsion. 
\end{abstract}

\section{Introduction}

In \cite{Ma1, Ma2} Yu.I. Manin introduced what is now called the Brauer--Manin obstruction.
To an element of the Brauer--Grothendieck group of a variety $X$ 
over a number field $k$ he attached a global reciprocity condition on the adelic points
of $X$ which is satisfied when an adelic point comes from a $k$-point.
In this paper we study the Brauer--Manin obstruction
on Kummer varieties, which are higher-dimensional generalisations of 
classical Kummer K3 surfaces.

Over complex numbers,
Kummer varieties in dimension greater than 2 were introduced in 1890 by W. Wirtinger \cite{W}.
Their topological and geometric properties were studied by A. Andreotti, E. Spanier
and K. Ueno, see \cite{Sp}, \cite{U1}, \cite{U2}.

Over non-closed fields, Kummer varieties come
not only from the quotients of abelian varieties by the antipodal involution, but also 
from the quotients of certain torsors. More precisely,
let $A$ be an abelian variety of dimension $g\geq 2$ over a field $k$ of characteristic
not equal to $2$. Let $Y$ be a $k$-torsor for $A$ whose class in
$\H^1(k,A)$ has order at most 2. Classically such torsors are referred to
as 2-{\em coverings} of $A$.
Kummer varieties considered in this paper are minimal desingularisations 
of the quotient $Y/\iota$ by the involution $\iota:Y\to Y$ induced by
the antipodal involution $[-1]:A\to A$. 
In the case $g=2$ we obtain Kummer
surfaces, a particular kind of K3 surface. Due to their intimate relation
to abelian varieties, Kummer surfaces are a popular testing
ground for conjectures on the geometry and arithmetic of K3 surfaces.
Rational points and Brauer groups of Kummer surfaces were studied in
\cite{SS, SZ2, IS, HS, CTTV, VAV}.

Rational points on Kummer varieties of higher dimension 
feature in the work of D. Holmes and R. Pannekoek \cite{HP}. Their result concerns
an abelian variety $A$ over a number field $k$: if
the set of $k$-points of the Kummer variety $X$
attached to $A^n$ is dense in the Brauer--Manin set of $X$, 
then there is a quadratic twist of $A$ over $k$ of rank at least $n$.
More recently, a Hasse principle for Kummer varieties, that are sufficiently general 
in an appropriate arithmetic sense, was established
conditionally on the finiteness of relevant Shafarevich--Tate groups 
by Y. Harpaz and one of the present authors \cite{HS}.
Somewhat surprisingly, the Brauer group does not show up in that statement.

Our aim in this paper is twofold. In Section \ref{geo} we establish geometric properties
of Kummer varieties analogous to similar properties of Kummer surfaces. 
We show, among other things, that the geometric Picard group $\Pic(\ov X)$ 
is a finitely generated free abelian group (Corollary \ref{cor0}).
In the characteristic zero case we describe a natural isomorphism of
Galois modules between the geometric Brauer group of a Kummer
variety and the geometric Brauer group of the corresponding abelian variety 
(Proposition \ref{Br}). From our previous result \cite{SZ1} we 
then deduce the finiteness of the quotient of
$\Br(X)$ by $\Br_0(X)=\Im[\Br(k)\to\Br(X)]$ when $k$ is finitely generated over $\Q$, 
see Corollary \ref{cor2}. Note, however, that the canonical class of a Kummer variety
of dimension $g\geq 3$ is represented by an effective divisor
(Proposition \ref{can}), thus higher-dimensional Kummer varieties
are not Calabi--Yau. Yonatan Harpaz asked if this could be relevant for
the tension which exists, in the light of the result of Holmes
and Pannekoek, between the heuristics for the ranks 
of elliptic curves over $\Q$ \cite{PPVW} and the conjecture that $\Q$-points of K3 surfaces
are dense in the Brauer--Manin set \cite[p. 77]{SkOb}, \cite[p. 484]{SZ1}.

The main goal of this paper is to study the Brauer group and the Brauer--Manin obstruction on
Kummer varieties. We prove the following general result.

\smallskip

\noindent{\bf Theorem \ref{g1}}
{\em Let $A$ be an abelian variety of dimension $>1$ over a number field $k$.
Let $X$ be the Kummer variety attached to a $2$-covering of $A$
such that $X(\A_k)\not=\emptyset$. Then $X(\A_k)^{\Br(X)_{\rm odd}}\not=\emptyset$,
where $\Br(X)_{\rm odd}\subset\Br(X)$ is the subgroup of elements of odd order.}

\smallskip

In Theorem \ref{t1} we give sufficient conditions on 
an abelian variety $A$ which guarantee
that the 2-torsion subgroup of $\Br(X)$ is contained in 
the algebraic Brauer group $\Br_1(X)=\Ker[\Br(X)\to\Br(\ov X)]$
and, moreover, $\Br_1(X)=\Br_0(X)$. 
The conditions of Theorem \ref{t1} are satisfied for
the Kummer variety $X$ attached to a $2$-covering of the Jacobian
of the hyperelliptic curve $y^2=f(x)$, where $f(x)\in k[x]$ is a
separable polynomial of odd degree $d\geq 5$ whose Galois group is the
symmetric or alternating group on $d$ letters. See Theorem \ref{t11}, where we
also treat products of Jacobians assuming that the splitting fields of the 
corresponding polynomials are linearly disjoint over $k$. 
This implies the following

\smallskip

\noindent{\bf Corollary \ref{g2}}
{\em Let $k$ be a number field.
Let $A$ be the product of Jacobians of elliptic or hyperelliptic curves 
$y^2=f_i(x)$, where $f_i(x)\in k[x]$ is a
separable polynomial of odd degree $m_i\geq 3$ whose Galois group is the 
symmetric group on $d$ letters. 
Assume that $\dim(A)>1$ and the splitting fields of the $f_i(x)$ are linearly disjoint over $k$.
If the Kummer variety $X$
attached to a $2$-covering of $A$ is everywhere locally soluble, then $X(\A_k)^\Br\not=\emptyset$.}

\smallskip

This explains the absence of the Brauer group from the statements of the Hasse principle
for K3 surfaces in Theorems A and B of \cite{HS}.

As a by-product of our calculations, we use a result of L. Dieulefait \cite{D} to
construct a Kummer K3 surface over $\Q$
of geometric Picard rank 17 with trivial Brauer group, see the examples at the end of the paper.
Previously known K3 surfaces with this property have geometric Picard rank 18, 19 and 20,
see \cite{SZ2, ISZ, IS}.

The first named author was partially supported by a grant EP/M020266/1
from EPSRC. The work on this paper started when he was
visiting the Pennsylvania State University and continued when he was at the Institute
for Advanced Study in Princeton, where he was supported by The Charles
Simonyi Endowment. He is grateful to both institutions. 
The second named author was partially supported by a grant
from the Simons Foundation (grant 246625 to Yuri Zarkhin).
A part of this work was done while he was a visitor at the Max-Planck-Institut
f\"ur Mathematik in Bonn, whose hospitality and support are gratefully acknowledged. We would like
to thank Yonatan Harpaz for extremely helpful comments and
Samir Siksek for directing us to the paper \cite{D}.
We are grateful to Tatiana Bandman and the referees for their comments.

\section{Kummer varieties and Kummer lattices} \label{geo}

Let $k$ be a field of characteristic different from 2 with an
algebraic closure $\bar k$  and the Galois group $\Ga=\Gal(\bar k/k)$.
For a variety $X$ over $k$ we write $\ov X=X\times_k\bar k$.
Let $A$ be an abelian variety over $k$ of dimension $g\geq 2$. 
We write $A^t$ for the dual abelian variety of $A$.

Let $T$ be a $k$-torsor for the group $k$-scheme $A[2]$.
We define the attached 2-covering of $A$ as the quotient 
$Y=(A\times_k T)/A[2]$ by the diagonal action of $A[2]$.
The first projection defines a morphism $f:Y\to A$
which is a torsor for $A[2]$ such that $T=f^{-1}(0)$.
The natural action of $A$ on $Y$ makes $Y$ a $k$-torsor for $A$.
In particular, there is an isomorphism of varieties $\ov Y\cong\ov A$.
Alternatively, $Y$ is the twisted form of $A$ defined by a 1-cocycle
with coefficients in $A[2]$ representing the class of $T$
in $\H^1(k,A[2])$, where $A[2]$ acts on $A$ by translations.

We have an exact sequence of $\Ga$-modules
\begin{equation}
0\lra A^t(\bar k)\lra\Pic(\ov Y)\lra\NS(\ov Y)\lra 0.\label{u2}
\end{equation}
The abelian groups $\NS(\ov Y)$ and $\NS(\ov A)$ are isomorphic. 
In fact, $\NS(\ov Y)$ and $\NS(\ov A)$ are also isomorphic
as $\Ga$-modules because
translations by the elements of $A(\bar k)$ act trivially on $\NS(\ov A)$, see \cite{Mu}.

The antipodal involution $\iota_A=[-1]:A\to A$ 
induces an involution $\iota_Y:Y\to Y$. 
It acts on $\Pic^0(\ov Y)=A^t(\bar k)$ as $[-1]$, which implies that
$$\H^0(\langle\iota_Y\rangle,A^t(\bar k))=A^t[2],\quad 
\H^1(\langle\iota_Y\rangle,A^t(\bar k))=0,$$
where we used the divisibility of $A^t(\bar k)$.
Taking the invariants of the action of 
$\iota_Y$ on the terms of (\ref{u2}) we obtain an exact sequence of $\Ga$-modules 
\begin{equation}
0\lra A^t[2]\lra \Pic(\ov Y)^{\iota_Y}\lra \NS(\ov Y)\lra 0. \label{b3}
\end{equation}

Let $\sigma:Y'\to Y$ be the blowing-up of the $2^{2g}$-point closed subscheme
$T\subset Y$. The involution $\iota_Y:Y\to Y$ preserves $T$ and so gives rise to an involution $\iota_{Y'}:Y'\to Y'$.

\bde
The {\bf Kummer variety} attached to $Y$ is the quotient $X=Y'/\iota_{Y'}$. 
\ede

By definition $\dim(X)=g\geq 2$. The fixed point set of $\iota_{Y'}$ is 
the exceptional divisor $E=\sigma^{-1}(T)$, which is 
a smooth divisor in $Y'$. A standard local calculation shows
that $X$ is smooth. Thus the natural surjective morphism $\pi:Y'\to X$ is a double covering 
whose branch locus is $E$. The divisor 
$\ov E=\sigma^{-1}(\ov T)$ is the disjoint union of $2^{2g}$ copies
of $\P^{g-1}_{\bar k}$. Let $D=\pi(E)\subset X$.

Let us now pause and describe some known facts about
Kummer varieties over $k=\C$. Spanier showed that these varieties
are simply connected \cite[Thm. 1]{Sp}. He also showed that their integral cohomology
groups are torsion-free, and computed the Betti numbers \cite[Thm. 2]{Sp}:
$$b_0=b_{2g}=1, \quad\quad b_{2i+1}=0, \quad \quad
b_{2i}=\left(\begin{array}{c}2g\\2i\end{array}\right)+2^{2g}, \ \text{where} \ \ 0< i< n.$$
The canonical class of $X$ was calculated by K. Ueno:
$$K_X=\frac{g-2}{2}[D],$$
see \cite[Lemma 16.11.1]{U2} or Proposition \ref{can} below. Since $K_X\geq 0$,
a theorem of K.~Ueno \cite[Prop. 3]{U1}, \cite[Thm. 16.2]{U2} says that
the Kodaira dimension of $X$ is $0$. 

We now return to the assumption that $k$ is an arbitrary field of characteristic 
different from 2.

\ble \label{lattice}
The subgroup of $\Pic(\ov X)$ generated by the classes
of the irreducible components of $\ov D$
is a free abelian group of rank $2^{2g}$ whose generators 
canonically correspond to the $\bar k$-points of $T$.
\ele
{\em Proof.} Let $E_i$, for $i=1,\ldots, 2^{2g}$, be the irreducible components
of $\ov E\subset \ov Y'$. Choose a line $L_i\cong\P^1_{\bar k}$ in each $E_i$.  
We define $D_i=\pi(E_i)\subset \ov X$, where $i=1,\ldots, 2^{2g}$.
The restriction of $\pi$ to $E_i$ is an isomorphism $E_i\to D_i$.

For $i\not=j$ we have $D_i\cap D_j=\emptyset$,
hence $([D_i].[\pi(L_j)])_{\ov X}=0$. 
The normal bundle $N$ to $E_i\cong\P^{g-1}_{\bar k}$ in $\ov Y'$ is $\O(-1)$.
By the standard formula
\cite[Prop. 2.6 (c)]{Fulton} for each $i=1,\ldots, 2^{2g}$ we have
$$([E_i].[L_i])_{\ov Y'}=(c_1(N).[L_i])_{E_i}=(\O(-1).[L_i])_{\P^{g-1}_{\bar k}}=-1.$$ 
Since $\pi^*[D_i]=2[E_i]$, by the projection formula we have
$$([D_i].[\pi(L_i)])_{\ov X}=(\pi^*[D_i].[L_i])_{\ov Y'}=-2.$$
Thus no non-trivial linear combination of the
classes $[D_i]$ is zero in $\Pic(\ov X)$. $\Box$

\medskip

We write $\Z[T]\subset \Pic(\ov X)$ for the subgroup 
described in Lemma \ref{lattice}. For $x\in T(\bar k)$ we denote the
corresponding generator of $\Z[T]$ by $e_x$.
Define $\Pi$ as the saturation of $\Z[T]$ in $\Pic(\ov X)$:
$$\Pi=\{x\in\Pic(\ov X)|\, nx\in\Z[T]\ \text{for some non-zero}\ n\in\Z\}.$$
For $g=2$ Nikulin proved in \cite[\S1]{N} that
$\Pi$ is a lattice in $\Q[T]=\Z[T]\otimes\Q$ generated
by $\Z[T]$ and the vectors $\frac{1}{2}\sum_{x\in H}e_x$,
where $H$ is a subset of $T(\bar k)\simeq A[2](\bar k)$ given by 
$L(x)=c$ for some $L\in\Hom(A[2],\F_2)$ and $c\in\F_2$. 
(This set of generators does not depend on the choice of an isomorphism 
$T(\bar k)\simeq A[2](\bar k)$ of $\bar k$-torsors for $A[2]$.)
We generalise this result to $g\geq 2$. In doing so we show that
$\Pic(\ov X)$ is torsion-free for any $g\geq 2$, see Proposition \ref{free} below.
In particular, $\Pi$ is also torsion-free, so $\Pi$ can be called the {\em Kummer lattice}.

\medskip

Write $Y_0=Y\setminus T$ and $X_0=\pi(\sigma^{-1}(Y_0))$.
Then $Y_0$ is the complement to a finite set in a smooth, proper and 
geometrically integral variety of dimension at least 2, so we have
\begin{equation}
\bar k[Y_0]=\bar k, \quad \Pic(\ov Y_0)=\Pic(\ov Y), \quad \Br(\ov Y_0)=\Br(\ov Y),
\label{u1} \end{equation}
where the last property follows from \cite[Cor. 6.2, p. 136]{Gr}. 

The involution $\iota_Y$ acts on $Y_0$ without fixed points,
hence $\pi:Y_0\to X_0=Y_0/\iota_Y$ is a torsor for $\Z/2$. There is a 
Hochschild--Serre spectral sequence \cite[Thm. III.2.20]{EC}
\begin{equation}
\H^p(\Z/2,\H^q_\et(\ov Y_0,\G_m))\Rightarrow\H^{p+q}_\et(\ov X_0,\G_m).
\label{ss}\end{equation}
Using (\ref{u1}) we deduce an exact sequence
\begin{equation}
0\lra \Z/2\lra \Pic(\ov X_0) \stackrel{\sigma_*\pi^*}\lra 
\Pic(\ov Y)^{\iota_Y}\lra 0, \label{b2}
\end{equation}
where the last zero is due to the fact that 
$\H^2(\Z/2,\bar k^*)=\bar k^*/\bar k^{*2}=0$ as $\char(k)\not=2$.
Using the fact that $\NS(\ov Y)\cong\NS(\ov A)$ is torsion-free,
we deduce from (\ref{b2}) and (\ref{b3}) 
a commutative diagram of $\Ga$-modules with exact rows and columns
\begin{equation}\begin{array}{ccccccccc}
&&&&0&&0&&\\
&&&&\uparrow&&\uparrow&&\\
&&&&\NS(\ov Y)&=&\NS(\ov Y)&&\\
&&&&\uparrow&&\uparrow&&\\
0&\lra& \Z/2&\lra& \Pic(\ov X_0) &\lra &\Pic(\ov Y)^{\iota_Y}&\lra& 0\\
&&||&&\uparrow&&\uparrow&&\\
0&\lra& \Z/2&\lra& \Pic(\ov X_0)_{\rm tors} &\lra& A^t[2]&\lra& 0\\
&&&&\uparrow&&\uparrow&&\\
&&&&0&&0&&
\end{array}\label{u3}
\end{equation}
Since $X$ is smooth, the natural restriction map $\Pic(\ov X)\to \Pic(\ov X_0)$
is surjective; thus $\Pic(\ov X_0)=\Pic(\ov X)/\Z[T]$.
This implies $\Pic(\ov X_0)_{\rm tors}=\Pi/\Z[T]$, so we obtain
a commutative diagram of $\Ga$-modules with exact rows and columns
\begin{equation}\begin{array}{ccccccccc}
&&&&0&&0&&\\
&&&&\uparrow&&\uparrow&&\\
&&&&\NS(\ov Y)&=&\NS(\ov Y)&&\\
&&&&\uparrow&&\uparrow&&\\
0&\lra& \Z[T]&\lra& \Pic(\ov X) &\lra &\Pic(\ov X_0)&\lra& 0\\
&&||&&\uparrow&&\uparrow&&\\
0&\lra& \Z[T]&\lra& \Pi &\lra& \Pic(\ov X_0)_{\rm tors}&\lra& 0\\
&&&&\uparrow&&\uparrow&&\\
&&&&0&&0&&
\end{array}\label{u4}
\end{equation}
For future reference we write the middle column of (\ref{u4}) as an exact sequence
of $\Ga$-modules
\begin{equation}
0\lra \Pi\lra\Pic(\ov X)\stackrel{\sigma_*\pi^*}\lra\NS(\ov Y)\lra 0. \label{a1}
\end{equation}

\bpr \label{free}
Let $X$ be a Kummer variety over a field of characteristic different from $2$.
Then the abelian group $\Pic(\ov X)$ is torsion-free. There is an isomorphism of
abelian groups $\Pic(\ov X_0)_{\rm tors}\cong A^t[2]\oplus\Z/2$.
\epr
{\em Proof.} The statements concern varieties over $\bar k$, so
we can assume that $X$ is attached to the trivial 2-covering $Y=A$.
The translations by points of order 2 commute with the antipodal
involution $[-1]:A\to A$. This implies that the finite commutative group 
$k$-scheme $\GG=A[2]\times_k\Z/2$ acts on $A$ so that the elements
of $A[2]$ act as translations and the generator of $\Z/2$ acts
as $[-1]$. It is easy to see that $\GG$ acts
freely on $A_1=A\setminus A[4]$ with quotient $A_1/\GG=X_0$. Hence the quotient morphism
$f:A_1\to X_0$ is a torsor for $\GG$. Since $g\geq 2$,
we have $\bar k[A_1]=\bar k$ and $\Pic(\ov A_1)=\Pic(\ov A)$.
The Cartier dual $\widehat\GG$ is isomorphic to $A^t[2]\times\Z/2$, 
so the exact sequence \cite[(2.5), p. 17]{Sk} gives an injective map 
$A^t[2]\oplus\Z/2\hookrightarrow\Pic(\ov X_0)$.
The bottom exact sequence of (\ref{u3}) shows that the cardinality
of $A^t[2]\oplus\Z/2$ equals the cardinality of $\Pic(\ov X_0)_{\rm tors}$, so we obtain an
isomorphism of abelian groups $A^t[2]\oplus\Z/2\tilde\lra\Pic(\ov X_0)_{\rm tors}$.

Since $\Z[T]$ is torsion-free, the natural map $\Pic(\ov X)\to\Pic(\ov X_0)$
induces an injective map of torsion subgroups. In particular, 
a non-zero torsion element of $\Pic(\ov X)$ is annihilated by 2 and corresponds
to a connected unramified double covering of $\ov X$. A double covering
of $\ov X$ is uniquely determined by its restriction to $\ov X_0$.
Therefore, it is enough to show that any connected unramified
double covering of $\ov X_0$ is a restriction of a {\em ramified} double covering of $\ov X$.
By the previous paragraph any such covering of $\ov X_0$ is of the form
$A_1/\HH$, where $\HH\subset\GG$ is a subgroup of index 2.

If $\HH=A[2]$, then $\ov A_1/A[2]=\ov A\setminus A[2]=\ov A_0$.
Write $\sigma:A'\to A$ for the blowing-up of $A[2]$ in $A$. Then the 
unramified double covering $\ov A_0\to\ov X_0$ extends to
the double covering $\ov A'\to\ov X$ ramified exactly in the
exceptional locus $\sigma^{-1}(A[2])$. 

If $\HH\not=A[2]$, then there is a non-zero $\phi\in\Hom(A[2],\Z/2)=A^t[2]$
such that $\HH$ is the
kernel of the homomorphism $A[2]\oplus\Z/2\to\Z/2$ given by
$(x,y)\mapsto \phi(x)$ or by $(x,y)\mapsto\phi(x)+y$.
Define $A_\phi=\ov A/\Ker(\phi)$. Choose $a\in A[2](\bar k)$ such that $\phi(a)\not=0$.
Then $\ov A_1/\HH$ is the quotient of
$A_\phi$ with $A_\phi[2]$ and $[2]^{-1}(\phi(a))$ removed,
by the involution $x\mapsto -x$ in the first case
and $x\mapsto\phi(a)-x$ in the second case. It follows that
the unramified double covering $\ov A_1/\HH\to \ov X_0$
is the restriction of the double covering of $\ov X$ ramified
in $\sigma^{-1}(A[2]\setminus \Ker(\phi))$ in the first case
and in $\sigma^{-1}(\Ker(\phi))$ in the second case. $\Box$

\bco \label{cor0}
Any Kummer variety $X$ of dimension $g\geq 2$
over a field $k$ of characteristic not equal to $2$
satisfies the following properties:

{\rm (i)} $\Pic^0(\ov X)=0$;

{\rm (ii)} $\Pic(\ov X)=\NS(\ov X)$ is torsion-free of rank $2^{2g}+\rk(\NS(\ov A))$;

{\rm (iii)} $\H^1_\et(\ov X,\Z_\ell)=0$ for any prime $\ell\not=\char(k)$;

{\rm (iv)} $\H^2_\et(\ov X,\Z_\ell)$ is torsion-free for any prime $\ell\not=\char(k)$.
\eco
{\em Proof.} Since $\Pic(\ov X)$ is torsion-free, we immediately obtain
(i) and $\Pic(\ov X)=\NS(\ov X)$. From diagram (\ref{u4}) we see that
the rank of this group is $2^{2g}+\rk(\NS(\ov Y))=2^{2g}+\rk(\NS(\ov A))$.
The Kummer sequence gives well-known isomorphisms
$$\H^1_\et(\ov X,\mu_{\ell^n})=\Pic(\ov X)[\ell^n]=0, \quad n\geq 1,$$
which imply (iii), by passing to the limit in $n$. The Kummer sequence also implies
the well-known fact that the torsion subgroup of $\H^2_\et(\ov X,\Z_\ell(1))$
coincides with the torsion subgroup of $\NS(\ov X)\otimes_\Z\Z_\ell$.
This gives (iv). $\Box$

\bco \label{cor1}
The Galois cohomology group $\H^1(k,\Pic(\ov X))$ is finite. 
The kernel of the natural map $\H^1(k,\Pic(\ov X))\to \H^1(k,\NS(\ov Y))$
is annihilated by $2$. If the order of the finite group $\H^1(k,\NS(\ov A))$
is a power of $2$, in particular, if $\NS(\ov A)$ is a trivial $\Ga$-module, then every 
element of odd order in $\Br_1(X)$ is contained in $\Br_0(X)$.
\eco
{\em Proof.} The finiteness of $\H^1(k,\Pic(\ov X))$ follows from 
the first statement of Proposition~\ref{free}.

The second statement of Proposition \ref{free} implies
that $\H^1(k,\Pic(\ov X_0)_{\rm tors})$ is annihilated by 2.
Since $\Z[T]$ is a permutation $\Ga$-module we have $\H^1(k,\Z[T])=0$.
By diagram (\ref{u4}) this implies that $\H^1(k,\Pi)$ is
a subgroup of $\H^1(k,\Pic(\ov X_0)_{\rm tors})$ and so is also annihilated by 2.
This proves the second statement.

Recall that $\NS(\ov Y)$ and $\NS(\ov A)$ are isomorphic
as $\Ga$-modules, so $\H^1(k,\NS(\ov Y))=\H^1(k,\NS(\ov A))$.
When the order of this group is a power of $2$, the order of $\H^1(k,\Pic(\ov X))$
is also a power of $2$.
The last statement is now immediate from the well known inclusion of $\Br_1(X)/\Br_0(X)$
into $\H^1(k,\Pic(\ov X))$. $\Box$

\medskip

We define $\Pi_1\subset \Pi$ as the kernel of the composed surjective map
$$\Pic(\ov X) \lra \Pic(\ov X_0)\lra \Pic(\ov Y)^{\iota_Y}.$$
Then $\Z[T]$ is a subgroup of $\Pi_1$ of index 2.
It is easy to see that $\Pi_1$
is generated by $\Z[T]$ and $\frac{1}{2}\sum_{x\in T(\bar k)}e_x$.
We thus have a canonical filtration
$$\Z[T]\subset \Pi_1\subset \Pi\subset \Pic(\ov X)$$
with successive factors $\Z/2$, $A^t[2]$, $\NS(\ov Y)=\NS(\ov A)$.
This filtration is respected by the action of $A[2]$ on $\ov Y$ and $\ov X$,
as well as by the action of the Galois group $\Ga$.

We summarise our discussion in the form of the following
commutative diagram with exact rows and columns,
where all arrows are group homomorphisms
which respect the actions of $\Ga$ and $A[2]$:
\begin{equation}\begin{array}{ccccccccc}
&&0&&0&&&&\\
&&\downarrow&&\downarrow&&&&\\
&&\Pi_1&=&\Pi_1&&&&\\
&&\downarrow&&\downarrow&&&&\\
0&\lra& \Pi&\lra&\Pic(\ov X)&\stackrel{\sigma_*\pi^*}\lra&\NS(\ov Y)&\lra& 0\\
&&\downarrow&&\downarrow&&||&&\\
0&\lra& A^t[2]&\lra&\Pic(\ov Y)^{\iota_Y}&\lra&\NS(\ov Y)&\lra& 0\\
&&\downarrow&&\downarrow&&\downarrow&&\\
&&0&&0&&0&&
\end{array} \label{diag}
\end{equation}

\noindent{\bf Remark 1} It is clear that $\Z[T]$ is a permutation $\Ga$-module.
Now consider the particular case when $T$ is a trivial torsor, i.e. $T\cong A[2]$
as $k$-torsors. The action of $\Ga$ on the set $A[2]$ fixes $0$.
It follows that not just $\Z[A[2]]$ but also $\Pi_1$ is a permutation $\Ga$-module.
Indeed, $\Pi_1$ has a $\Ga$-stable $\Z$-basis consisting of $e_x$ for 
$x\in A[2]\setminus\{0\}$ and $\frac{1}{2}\sum_{x\in A[2]}e_x$.
Note, however, that this basis is not $A[2]$-stable.

\medskip

The following proposition shows that the canonical class of 
a Kummer variety of dimension $g\geq 3$
is represented by an effective divisor, so such varieties are not Calabi--Yau.
In the case ${\rm char}(k)=0$, this was proved in \cite[Lemma 16.11.1]{U2}.

\bpr \label{can}
We have $K_{\ov X}=\frac{g-2}{2}[\ov D]=\frac{g-2}{2}\sum_{x\in T(\bar k)}e_x$.
\epr
{\em Proof.} The natural map $\pi^*:\Pic(\ov X)\to\Pic(\ov Y')$ is injective.
Indeed, its kernel is contained in $\Pi_1$, by the exactness of the middle column of (\ref{diag}).
In the notation of the proof of Lemma \ref{lattice} we have $\pi^*[D_i]=2[E_i]$,
hence $\pi^*$ is injective on $\Pi_1$. Since $K_{\ov Y}=0$, the standard formulae give
$K_{\ov Y'}=(g-1)\sum [E_i]$ and $K_{\ov Y'}=\pi^*K_{\ov X}+\sum [E_i]$.
Now our statement follows from
the injectivity of $\pi^*:\Pic(\ov X)\to\Pic(\ov Y')$. $\Box$

\bpr \label{Br}
Assume that the characteristic of $k$ is zero.
Then the morphisms $\pi:Y'\to X$ and $\sigma:Y'\to Y$ induce isomorphisms of
$\Ga$-modules
$$\Br(\ov X)\tilde\lra\Br(\ov Y')\tilde\longleftarrow \Br(\ov Y)\cong\Br(\ov A).$$
\epr
{\em Proof} The last isomorphism is due to the fact
that $Y$ is the twist of $A$ by a 1-cocycle with coefficients
in $A[2]$, but the induced action of $A[2]$ on $\Br(\ov A)$ is trivial. 
In fact, the whole group $A(\bar k)$ acts trivially on the finite
group $\Br(\ov A)[n]$ for every integer $n$, because any
homomorphism from the divisible group $A(\bar k)$ to the finite
group $\Aut(\Br(\ov A)[n])$ is trivial.

The middle isomorphism is a consequence of the birational invariance of
the Brauer group of a smooth and projective variety over a field
of characteristic zero. 

The natural map $\pi^*:\Br(\ov X)\to\Br(\ov Y')$ is a map of $\Ga$-modules.
To prove that it is an isomorphism we can work 
over an algebraically closed field of characteristic zero
and so assume that $Y=A$. We remark that
Grothendieck's exact sequence \cite[Cor. 6.2, p. 137]{Gr} gives
an exact sequence
$$0\lra\Br(\ov X)\lra\Br(\ov X_0)\lra \bigoplus\H^1(\P^{g-1}_{\bar k},\Q/\Z),$$
where the terms in the direct sum are numbered by the $2^{2g}$ points of
$A[2](\bar k)$. We have $\H^1(\P^{g-1}_{\bar k},\Z/n)=0$ for any
positive integer $n$, so the natural map $\Br(\ov X)\tilde\lra\Br(\ov X_0)$
is an isomorphism. By (\ref{u1}) there is a natural isomorphism 
$\Br(\ov A)\tilde\lra\Br(\ov A_0)$.

We analyse the map $\pi^*:\Br(\ov X_0)\to\Br(\ov A_0)$ using the 
spectral sequence (\ref{ss}). We have already seen that
$\H^2(\Z/2,\bar k^*)=0$. We have a natural isomorphism
$\Pic(\ov A_0)=\Pic(\ov A)$ and we claim that $\H^1(\Z/2,\Pic(\ov A))=0$.
In view of the exact sequence (\ref{u2}) it is enough to prove that 
$\H^1(\Z/2,A^t)=\H^1(\Z/2,\NS(\ov A))=0$.
The torsion-free group $\NS(\ov A)$ is a subgroup of 
$\H^2_\et(\ov A,\Z_\ell(1))$ for any prime $\ell$.
The involution $[-1]$ acts trivially on $\H^2_\et(\ov A,\Z_\ell(1))$
and hence on $\NS(\ov A)$. It follows that $\H^1(\Z/2,\NS(\ov A))=0$.
On the other hand, $[-1]$ acts on $\Pic^0(\ov A)\cong A^t(\bar k)$ as $[-1]$, implying
$\H^1(\Z/2,A^t)=0$. The spectral sequence (\ref{ss}) now gives
an injective map $\Br(\ov X)\hookrightarrow\Br(\ov A)$.

By the well known Grothendieck's computation \cite[\S 8]{Gr}
we have $\Br(\ov A)\cong(\Q/\Z)^{b_2-\rho}$, where 
$b_2=g(2g-1)$ is the dimension of $\H^2_\et(\ov A,\Q_\ell(1))$ and $\rho=\rk(\NS(\ov A))$. To complete the proof it is enough to show that 
the corank of the divisible part of $\Br(\ov X)$ is $g(2g-1)-\rho$. 
(Indeed, any injective homomorphism 
$(\Q/\Z)^r\to (\Q/\Z)^r$ is an isomorphism.)
By Corollary \ref{cor0} (ii) $\Pic(\ov X)$ is torsion-free of rank $\rho+2^{2g}$.
Thus it remains to show that the dimension of
$\H^2_\et(\ov X,\Q_\ell(1))$ is $g(2g-1)+2^{2g}$
for any prime $\ell$. The Gysin sequence gives an exact sequence
$$0\lra (\Q_\ell)^{2^{2g}}\lra \H^2_\et(\ov X,\Q_\ell(1))\lra 
\H^2_\et(\ov X_0,\Q_\ell(1))\lra 0.$$
The spectral sequence $\H^p(\Z/2,\H^q_\et(\ov A_0,\Q_\ell))
\Rightarrow\H^{p+q}_\et(\ov X_0,\Q_\ell)$ degenerates because
each $\H^q_\et(\ov A_0,\Q_\ell)$ is a vector space over a field of characteristic $0$.
We obtain
$$\H^n_\et(\ov X_0,\Q_\ell)=\H^n_\et(\ov A_0,\Q_\ell)^{[-1]^*}$$
for all $n\geq 0$. In particular, the dimension of 
$\H^2_\et(\ov X_0,\Q_\ell(1))$ is $g(2g-1)$, as required. $\Box$

\bco \label{cor2}
Let $k$ be a field finitely generated over $\Q$. Let $X$ be the Kummer
variety attached to a $2$-covering of an abelian variety. Then
the groups $\Br(X)/\Br_0(X)$ and $\Br(\ov X)^\Ga$ are finite.
\eco
{\em Proof.} By
the spectral sequence $\H^p(k,\H^q_\et(\ov X,\G_m))
\Rightarrow \H^{p+q}_\et(X,\G_m)$ and Corollary
\ref{cor1} the finiteness of $\Br(\ov X)^\Ga$ implies
the finiteness of $\Br(X)/\Br_0(X)$. By Proposition \ref{Br}
this follows from the finiteness of $\Br(\ov A)^\Ga$ which is
established in \cite{SZ1}. $\Box$

\medskip

\noindent{\bf Remark 2} Assume that ${\rm char}(k)=0$.
The commutative diagram
$$\begin{array}{ccc} \Br(\ov X)&\tilde\lra&\Br(\ov Y)\\
\uparrow&&\uparrow\\
\Br(X)&\stackrel{\sigma_*\pi^*}\lra&\Br(Y)\end{array}$$
identifies $\Br(X)/\Br_1(X)$ with a subgroup of $\Br(Y)/\Br_1(Y)$.

\section{When the Hasse principle is unobstructed}

Let $n$ be an {\em odd} integer and let $k$ be a field of characteristic
coprime to $2n$. 
Let $\Lambda$ be a $\Ga$-module such that $n\Lambda=0$.

If $A$ be an abelian variety over $k$, then
$[-1]$ acts on $\H^q_\et(\ov A,\Lambda)$ by $(-1)^q$,
where $q\geq 0$. Hence for a $2$-covering $Y$ of $A$ the involution
$\iota_Y$ acts on $\H^q_\et(\ov Y,\Lambda)$ by $(-1)^q$.

For $m\geq 0$ let $\H^m_\et(Y,\Lambda)^+$ be the $\iota_Y$-invariant 
subgroup of $\H^m_\et(Y,\Lambda)$.
Let $\H^m_\et(Y,\Lambda)^-$ be the $\iota_Y$-anti-invariant subgroup,
i.e. the group of elements on which $\iota_Y$ acts by $-1$.
Since $n$ is odd, we can write
\begin{equation}
\H^m_\et(Y,\Lambda)=\H^m_\et(Y,\Lambda)^+\oplus \H^m_\et(Y,\Lambda)^-.
\label{dec}
\end{equation}

\bpr \label{coh}
Let $Y$ be a $2$-covering of an abelian variety $A$.
Then we have a canonical decomposition of abelian groups
$$\H^2_\et(Y,\Lambda)=\H^2(k,\Lambda)\oplus \H^1(k,\H^1_\et(\ov Y,\Lambda))
\oplus \H^2_\et(\ov Y,\Lambda)^\Ga$$
compatible with the natural action of 
the involution $\iota_Y$ on $\H^2_\et(Y,\Lambda)$, so that 
$$\H^2_\et(Y,\Lambda)^+=\H^2(k,\Lambda)\oplus \H^2_\et(\ov Y,\Lambda)^\Ga \quad
\text{and}\quad\H^2_\et(Y,\Lambda)^-=\H^1(k,\H^1(\ov Y,\Lambda)).$$
\epr
{\em Proof.} Let $m\geq 1$. The morphisms $Y\to\Spec(k)$ and $\ov Y\to Y$
induce the $\iota_Y$-equivariant maps 
$$\alpha_m:\H^m(k,\Lambda)\lra \H^m_\et(Y,\Lambda), \quad
\beta_m:\H^m_\et(Y,\Lambda)\lra \H^m_\et(\ov Y,\Lambda)^\Ga,
\quad \beta_m\alpha_m=0.$$ 
Since $\iota_Y$ acts trivially on $\H^m(k,\Lambda)$
and on $\H^2_\et(\ov Y,\Lambda)$, we have
$$\Im(\alpha_m)\subset \H^m_\et(Y,\Lambda)^+,\quad\quad
\H^2_\et(Y,\Lambda)^-\subset\Ker(\beta_2).$$ 
We claim that it is enough to show that

\smallskip

$\alpha_m:\H^m(k,\Lambda)\to \H^m_\et(Y,\Lambda)^+$ 
has a retraction, for every $m\geq 0$;

\smallskip

$\beta_2:\H^2_\et(Y,\Lambda)^+\to \H^2_\et(\ov Y,\Lambda)^\Ga$ 
has a section.

\smallskip

\noindent Indeed, if this is true, then $\H^m(k,\Lambda)$ is a direct summand of 
$\H^m_\et(Y,\Lambda)^+$ for $m\geq 0$. Moreover, $\H^2_\et(\ov Y,\Lambda)^\Ga$
is a direct summand of $\H^2_\et(Y,\Lambda)^+$, so that
$$\H^2_\et(Y,\Lambda)=\H^2(k,\Lambda)\oplus \Ker(\beta_2)/\Im(\alpha_2)\oplus
\H^2_\et(\ov Y,\Lambda)^\Ga.$$
Note that the maps $\alpha_m$ and $\beta_m$ are the canonical maps in the spectral sequence
\begin{equation}
\H^p(k,\H^q_\et(\ov Y,\Lambda))\Rightarrow\H^{p+q}_\et(Y,\Lambda). \label{leray}
\end{equation}
From (\ref{leray}) we obtain the exact sequence 
$$0\lra \Ker(\beta_2)/\Im(\alpha_2)\lra \H^1(k,\H^1(\ov Y,\Lambda))\lra \Ker(\alpha_3)=0.$$
Since $\iota_Y$ acts on $\H^1(\ov Y,\Lambda)$ by $-1$ and
the Galois group $\Ga$ commutes with $\iota_Y$, we get
$$\Ker(\beta_2)/\Im(\alpha_2)=\H^1(k,\H^1(\ov Y,\Lambda))\subset \H^2_\et(Y,\Lambda)^-.$$
Now the claim follows from (\ref{dec}).

Let us construct a retraction of $\alpha_m$.
Recall that $T$ is a 0-dimensional subscheme of $Y$, and so we have
a restriction map $\H^m_\et(Y,\Lambda)\to \H^m_\et(T,\Lambda)$. Write
$T$ as a disjoint union of closed points
$$T=\bigsqcup_{i=1}^r\Spec(k_i),$$
where each $k_i$ is a finite field extension of $k$.
Then $\H^m_\et(T,\Lambda)$ is the direct sum of the Galois cohomology
groups $\H^m(k_i,\Lambda)$ for $i=1,\ldots,r$.
The composition of the restriction $\H^m(k,\Lambda)\to \H^m(k_i,\Lambda)$
and the corestriction $\H^m(k_i,\Lambda)\to \H^m(k,\Lambda)$ is the multiplication by $[k_i:k]$.
The direct sum of these corestriction maps is a map
$\H^m_\et(T,\Lambda)\to \H^m(k,\Lambda)$ whose composition with the natural
restriction map $\H^m(k,\Lambda)\to \H_\et^m(T,\Lambda)$ is the multiplication
by $|T(\bar k)|=2^{2g}$. Since $n$ is odd, there is an integer $r$ such that
$2^{2g}r\equiv 1 \bmod n$.
Thus the composition 
\begin{equation}
\H^m_\et(Y,\Lambda)\lra \H^m_\et(T,\Lambda)\lra
\H^m(k,\Lambda)\stackrel{[r]}\lra\H^m(k,\Lambda)\label{retract}
\end{equation}
is a retraction of $\alpha_m$. 
Since $T\subset Y$ is the fixed point set
of $\iota_Y$, this retraction is $\iota_Y$-equivariant.

Let us construct a section of $\beta_2$. 
The translations by the points of $A(\bar k)$ act trivially on $\H^q_\et(\ov A,\Lambda)$
for any $q\geq 0$, so
we have canonical isomorphisms of $\Ga$-modules
\begin{equation}
\H^2_\et(\ov Y,\Lambda)=\H^2_\et(\ov A,\Lambda)=\Hom(\wedge^2 A[n],\Lambda).
\label{g7}
\end{equation}
Since $[Y]\in\H^1(k,A)[2]$ and $n$ is odd, the multiplication by $n$ on $A$
defines a morphism $[n]:Y\to Y$ which is a torsor with structure group $A[n]$.
We denote this torsor by $\T_n$. The class of this torsor is an element $[\T_n]\in \H^1_\et(Y,A[n])$.
Using cup-product we obtain a class
$$\wedge^2[\T_n]\in \H^2_\et(Y,\wedge^2 A[n]).$$
The isomorphisms (\ref{g7}) give rise to a pairing
$$\H^2_\et(Y,\wedge^2 A[n])\times \H^2_\et(\ov Y,\Lambda)^\Ga\lra \H^2_\et(Y,\Lambda).$$
Let $s:\H^2_\et(\ov Y,\Lambda)^\Ga\to \H^2_\et(Y,\Lambda)$ be the map
defined by pairing with the class $\wedge^2[\T_n]$. As $n$ is odd,
the same proof as in \cite[Prop. 2.2]{SZ2} (where we treated the case $Y=A$)
shows that $s$ is a section of the natural map
$\H^2_\et(Y,\Lambda)\to \H^2_\et(\ov Y,\Lambda)^\Ga$. 

Since $[-1]\cdot[n]=[-n]=[n]\cdot[-1]$ we see that $\iota_Y^*[\T_n]=[\T_{-n}]$,
and the torsor $\T_{-n}$ is obtained from $\T_n$ by applying the automorphism
$[-1]$ to the structure group $A[n]$. Hence $\iota_Y^*[\T_n]=-[\T_n]$.
It follows that $\iota_Y^*(\wedge^2[\T_n])=\wedge^2[\T_n]$. We conclude
that $s$ is $\iota_Y$-equivariant, and so is a section of 
$\beta_2:\H^2_\et(Y,\Lambda)^+\to \H^2_\et(\ov Y,\Lambda)^\Ga$. $\Box$

\medskip

\noindent{\bf Remark 3} The restriction of the morphism $[n]:Y\to Y$ to $T$
has a section given by the identity map $T\tilde\lra T$. Thus the restriction of $\T_n$
to $T$ is trivial. It follows that $s(\H^2_\et(\ov Y,\Lambda)^\Ga)$
is contained in the kernel of the restriction map
$\H^2_\et(Y,\Lambda)\to \H^2_\et(T,\Lambda)$.

\medskip

Recall that the Kummer variety $X$ attached to $Y$ is defined as follows.
Let $\sigma:Y'\to Y$ be the blowing-up of $T$ in $Y$. Then $\pi:Y'\to X$ is the
double cover which is the quotient by a natural involution on $Y'$
compatible with $\iota_Y$. We note that the same variety
$X$ can also be obtained from any quadratic twist of $Y$. More precisely, let 
$F$ be an \'etale $k$-algebra of dimension 2, i.e. $k\oplus k$ or a quadratic extension of $k$.
We denote by $A_F$ the quadratic twist of the abelian variety $A$ by $F$, defined as the quotient
of $A\times_k\Spec(F)$ by the simultaneous action of $\Z/2$ such that the generator of $\Z/2$
acts on $A$ as $[-1]$ and on $\Spec(F)$ as $c\in\Gal(F/k)$, $c\not=0$.
(In the case of $F=k\oplus k$ the action of $c$ permutes the factors of $\Spec(F)$,
so that $A_F=A$ in this case.)
We define $Y_F$ similarly, replacing $[-1]$ with $\iota_Y$. 
Since $[-1]$ commutes with translations by the elements of $A[2]$, 
we have a morphism $Y_F\to Y_F/A[2]=A_F$, which is a 2-covering of $A_F$ defined by the same
$k$-torsor $T$ for $A[2]=A_F[2]$. 
The blowing-up $\sigma_F:Y'_F\to Y_F$ of the closed subscheme $T\subset Y_F$
has an involution compatible with $\iota_{Y_F}$. It 
gives rise to the double covering $\pi_F:Y'_F\to X$.

Let ${Y_F}_ 0=Y_F\setminus T$.
We have a commutative diagram
\begin{equation}\begin{array}{ccccc}
\H^2_\et(X,\Lambda)&\stackrel{\pi_F^*}\lra &
\H^2_\et(Y'_F,\Lambda)&\stackrel{\sigma_{F*}}\lra& \H^2_\et(Y_F,\Lambda)\\
\downarrow&&\downarrow&&||\\
\H^2_\et(X_0,\Lambda)&\stackrel{\pi_F^*}\lra &
\H^2_\et({Y_F}_0,\Lambda)&=& \H^2_\et({Y_F}_0,\Lambda)\end{array}
\label{h}
\end{equation}
The restriction map
$\H^2_\et(Y_F,\Lambda)\tilde\lra\H^2_\et({Y_F}_ 0,\Lambda)$
is an isomorphism 
by the purity of \'etale cohomology \cite[Remark VI.5.4 (b)]{EC}
as $\codim_Y (T)\geq 2$.
The map $\sigma_{F*}$ is the composition of the restriction to 
the open set ${Y_F}_ 0\subset Y'_F$
and the inverse of $\H^2_\et(Y_F,\Lambda)\tilde\lra\H^2_\et({Y_F}_ 0,\Lambda)$.
In particular, the composition 
$$\H^2_\et(Y_F,\Lambda)\stackrel{\sigma_F^*}\lra \H^2_\et(Y'_F,\Lambda)
\stackrel{\sigma_{F*}}\lra \H^2_\et(Y_F,\Lambda)$$
is the identity map.

For the sake of completeness we note that the
Hochschild--Serre spectral sequence \cite[Thm. III.2.20]{EC}
$$\H^p(\Z/2,\H^q_\et({Y_F}_0,\Lambda))\Rightarrow \H^{p+q}_\et(X_0,\Lambda)$$
gives canonical isomorphisms
$$\H^p_\et(X_0,\Lambda)\tilde\lra\H^p_\et({Y_F}_0,\Lambda)^+.$$
Indeed, $\H^p(\Z/2,\H^q_\et({Y_F}_0,\Lambda))=0$
for $p\geq 1$, since $2$ and $n$ are coprime.

\bpr \label{propo}
Let $X$ be the Kummer variety attached to
a $2$-covering $Y$ of an abelian variety of dimension at least $2$.
Let $n\geq 1$ be an odd integer.
For any $x\in \H^2_\et(X,\mu_n)$ there exists an $a_0\in\H^2(k,\mu_n)$ 
such that for any \'etale $k$-algebra $F$ of dimension $2$ we have
$$\sigma_{F*}\pi_F^*(x)-a_0\in s(\H^2_\et(\ov Y_F,\mu_n)^\Ga),$$
where $s$ is the section of the natural map 
$\H^2_\et(Y_F,\mu_n)\to\H^2_\et(\ov Y_F,\mu_n)^\Ga$ constructed 
in the proof of Proposition \ref{coh}.
\epr
{\em Proof.} Since $\sigma_{F*}\pi_F^*(x)$ is $\iota_{Y_F}$-invariant,
by Proposition \ref{coh} we have $\sigma_{F*}\pi_F^*(x)=a_0+s(a)$
for some $a_0\in\H^2(k,\mu_n)$ and $a\in \H^2_\et(\ov Y_F,\mu_n)^\Ga$. 
We need to show that $a_0$ does not depend on $F$.
Recall that 
$$X\setminus X_0\cong Y'_F\setminus {Y_F}_0=
\P^{g-1}_T=\P^{g-1}_k\times_k T,$$
and the natural morphism $Y'_F\setminus {Y_F}_0\to Y_F\setminus {Y_F}_0=T$
is the structure morphism $\P^{g-1}_T\to T$. 
We have a commutative diagram, where the vertical arrows are the natural
restriction maps
\begin{equation}
\xymatrix{\H^2_\et(X,\mu_n) \ar[r]^{\pi_F^*}\ar[d]^{\tau_1}&
\H^2_\et(Y'_F,\mu_n)\ar[d]^{\tau_2}&
\H^2_\et(Y_F,\mu_n)\ar[l]_{\sigma_F^*}\ar[d]^{\tau_3}\\
\H^2_\et(\P^{g-1}_T,\mu_n) \ar[r]^{\rm id}&\H^2_\et(\P^{g-1}_T,\mu_n)
\ar@<-.5ex>[r]_{\rho^*} &
\H^2_\et(T,\mu_n) \ar@<-.5ex>[l]}\label{h1}
\end{equation}
A choice of a $k$-point in $\P^{g-1}_k$ defines a 
section $\rho$ of the structure morphism $\P^{g-1}_T\to T$,
and we denote by $\rho^*$ the induced map
$\H^2_\et(\P^{g-1}_T,\mu_n)\to \H^2_\et(T,\mu_n)$.

Recall from (\ref{retract}) that $a_0$ is obtained by
applying to $\tau_3\sigma_{F*}\pi_F^*(x)$ the corestriction map 
$\H^2_\et(T,\mu_n)\to\H^2(k,\mu_n)$ followed by the multiplication by $r$.

Let $y=\sigma_F^*\sigma_{F*}\pi_F^*(x)-\pi_F^*(x)\in \H^2_\et(Y'_F,\mu_n)$.
We claim that for any closed point
$i:\Spec(K)\hookrightarrow Y'_F$ we have 
$$i^*(y)=0\in \H^2(K,\mu_n).$$
Indeed, $\sigma_{F*}\sigma_F^*={\rm id}$ implies that
$\sigma_{F*}(y)=0$ and hence $y$ goes to 0
under the restriction map $\H^2_\et(Y'_F,\mu_n)\to 
\H^2_\et({Y_F}_ 0,\mu_n)$. The natural injective map of \'etale sheaves
$\mu_n\to\G_m$ gives rise to the canonical maps
$\H^2_\et(Y'_F,\mu_n)\to\Br(Y'_F)$ and 
$\H^2_\et({Y_F}_ 0,\mu_n)\to \Br({Y_F}_ 0)$.
By Grothendieck's purity theorem for the Brauer group \cite[III, \S 6]{Gr}
the natural restriction map $\Br(Y'_F)\to \Br({Y_F}_ 0)$
is injective. Hence the image of $y$ in $\Br(Y'_F)$ is zero.
On the other hand, the map $\H^2(K,\mu_n)\to \Br(K)$
is injective by Hilbert's Theorem 90. This implies $i^*(y)=0$.

In particular, we have $\rho^*\tau_2(y)=0$, hence
$\rho^*\tau_2\pi_F^*(x)=\rho^*\tau_2\sigma_F^*\sigma_{F*}\pi_F^*(x)$.
The commutativity of the right hand square of (\ref{h1}) and
the fact that $\rho$ is a section
of the structure morphism $\P^{g-1}_T\to T$ imply
that $\rho^*\tau_2\sigma_F^*=\tau_3$. Hence 
$\rho^*\tau_2\pi_F^*(x)=\tau_3\sigma_{F*}\pi_F^*(x)$.
By the commutativity of the left hand square of (\ref{h1}) 
we have $\tau_2\pi_F^*(x)=\tau_1(x)$.
Hence $\rho^*\tau_2\pi_F^*(x)=\rho^*\tau_1(x)$, which does not depend on $F$.
We conclude that $\tau_3\sigma_{F*}\pi_F^*(x)$,
and hence also $a_0$, do not depend on $F$. $\Box$

\bigskip

Now let $k$ be a number field. We write $\A_k$ for the ring of ad\`eles of $k$.
If $X$ is a proper variety over $k$ we have $X(\A_k)=\prod X(k_v)$,
where $v$ ranges over all places of $k$. The Brauer--Manin
pairing $X(\A_k)\times\Br(X)\to\Q/\Z$ is given by the sum of local invariants
of class field theory, see \cite[\S 5.2]{Sk}.
For a subgroup $B\subset\Br(X)$ we denote by $X(\A_k)^B\subset
X(\A_k)$ the orthogonal complement to $B$ under this pairing.

\bthe \label{g1}
Let $A$ be an abelian variety of dimension $g\geq 2$ over a number field $k$.
Let $X$ be the Kummer variety attached to a $2$-covering of $A$
such that $X(\A_k)\not=\emptyset$. Then $X(\A_k)^{\Br(X)_{\rm odd}}\not=\emptyset$, where $\Br(X)_{\rm odd}\subset\Br(X)$ 
is the subgroup of elements of odd order.
\ethe
{\em Proof.} By Corollary \ref{cor2} the group $\Br(X)/\Br_0(X)$ is finite.
It follows that $\Br(X)_{\rm odd}$ is generated by finitely many elements
modulo $\Br(X)_{\rm odd}\cap\Br_0(X)$. Hence there is an odd integer $n$ such that 
the images of $\Br(X)_{\rm odd}$ and $\Br(X)[n]$ in $\Br(X)/\Br_0(X)$ are equal.
Since the sum of local invariants of an element of $\Br_0(X)$ is always zero, this implies
that $X(\A_k)^{\Br(X)_{\rm odd}}=X(\A_k)^{\Br(X)[n]}$.

We have the natural maps
$$\H^2_\et(X,\G_m)\stackrel{\pi^*}\lra 
\H^2_\et(Y',\G_m)\stackrel{\sigma_*}\lra \H^2_\et(Y,\G_m).$$
Here $\sigma_*$ is the composition of the restriction to 
the open set $Y_0\subset Y'$ and the inverse of the restriction
$\H^2_\et(Y,\G_m)\tilde\lra\H^2_\et(Y_0,\G_m)$, which is an isomorphism by
Grothendieck's purity theorem for the Brauer group \cite[III, Cor. 6.2]{Gr} as
$\codim_Y (T)\geq 2$. These maps are compatible with
the similar maps (\ref{h}) with finite coefficients $\Lambda=\mu_n$.
Now the Kummer sequences for $X$ and $Y$ give rise to the commutative diagram
$$\xymatrix{
\H^2_\et(Y,\mu_n)\ar[r]&\Br(Y)[n]\ar[r]&0\\
\H^2_\et(X,\mu_n)\ar[r]\ar[u]^{\sigma_*\pi^*}&\Br(X)[n]\ar[r]\ar[u]^{\sigma_*\pi^*}&0
}$$
The same considerations apply if we replace $Y$ by any quadratic twist $Y_F$.

Take any $\sA\in\Br(X)[n]$ and lift it to some $x\in \H^2_\et(X,\mu_n)$.
By the commutativity of the previous diagram
$\sigma_{F*}\pi_F^*(\sA)\in \Br(Y_F)[n]$ comes from 
$\sigma_{F*}\pi_F^*(x)\in \H^2_\et(Y_F,\mu_n)$. By
Proposition \ref{propo} there is $a_0\in\H^2(k,\mu_n)$ such that
$\sigma_{F*}\pi_F^*(x)-a_0\in s(\H^2_\et(\ov Y_F,\mu_n)^\Ga)$.

Now we can complete the proof of the theorem. Let $(P_v)\in X(\A_k)$.
For each $v$ there is a class $\alpha_v\in\H^1(k_v,\mu_2)=k_v^*/k_v^{*2}$
such that $P_v$ lifts to a $k_v$-point on the quadratic twist $Y'_{k_v(\sqrt{\alpha_v})}$,
which is a variety defined over $k_v$.
By weak approximation in $k$ we can assume that $\alpha_v$ comes from
$\H^1(k,\mu_2)=k^*/k^{*2}$, and hence assume that 
$Y'_{k_v(\sqrt{\alpha_v})}\cong Y'_F\times_k k_v$ for some 
\'etale $k$-algebra $F$ of dimension 2. 

It follows that $Y'_F(k_v)\not=\emptyset$ and hence $Y_F(k_v)\not=\emptyset$. 
Since $Y_F$ is smooth, the non-empty set $Y_F(k_v)$ is Zariski dense in $Y_F$.
Thus there is a $k_v$-point $R_v\in Y_F(k_v)$ such that the point
$M_v=[n]R_v\in Y_F(k_v)$ is contained in the open subset ${Y_F}_0$. 
The specialisation of $\T_n$ at $M_v$ contains a $k_v$-point, hence is a trivial torsor.
It follows that the specialisation of the class $\wedge^2[\T_n]\in \H^2(Y_F,\wedge^2 A[n])$ 
at $M_v$ is zero. By the construction of the section $s$ in the proof of Proposition \ref{coh}
we obtain that $s(a)\in \H^2(Y_F,\mu_n)$ evaluated at $M_v$ is zero for any $a$.
Therefore, $\sigma_{F*}\pi_F^*(x)$ evaluated at $M_v$ is the image of
$a_0$ in $\H^2(k_v,\mu_n)$, and hence $\sigma_{F*}\pi_F^*(\sA)(M_v)\in \Br(k_v)$
comes from a global element $a_0\in \H^2(k,\mu_n)=\Br(k)[n]$.

Since $M_v\in {Y_F}_0(k_v)$ there exists a unique point
$M'_v\in Y'_F(k_v)$ such that $\sigma_F(M'_v)=M_v$.
Let $Q_v=\pi_F(M'_v)\in X(k_v)$.
By the projection formula we have $\sA(Q_v)=\sigma_{F*}\pi_F^*(\sA)(M_v)$. Since
this is the image of $a_0\in \Br(k)$ under the restriction map to $\Br(k_v)$,
the sum of local invariants of $\sA$ evaluated at the adelic point
$(Q_v)\in X(\A_k)$ is zero. 
Thus $(Q_v)\in X(\A_k)^{\Br(X)[n]}=X(\A_k)^{\Br(X)_{\rm odd}}$. $\Box$

\section{Kummer varieties attached to products of abelian varieties} \label{kuku}

For an abelian group $G$ we denote by $G\{\ell\}$ the $\ell$-primary subgroup of $G$.

\bpr \label{duo}
Let $k$ be a field and let $\ell$ be a prime different from the
characteristic of $k$. 
Let $A_1,\ldots,A_n$ be principally polarised abelian varieties over $k$ 
such that the fields $k(A_i[\ell])$ are pairwise linearly disjoint over $k$, 
where $i=1,\ldots,n$. Assume that each $\Ga$-module $A_i[\ell]$
is simple, and, moreover, if $\dim(A_i)>1$, then it is absolutely simple.
For any $2$-covering $Y$ of $A=\prod_{i=1}^n A_i$ 
we have $\Br(Y)\{\ell\}\subset\Br_1(Y)$, in particular, 
$\Br(A)\{\ell\}\subset\Br_1(A)$. If ${\rm char}(k)=0$ and $\dim(A)\geq 2$, 
for the Kummer variety $X$ attached to $Y$ we have $\Br(X)\{\ell\}\subset\Br_1(X)$.
\epr
{\em Proof.} Let $m$ be a positive integer.
The Kummer sequences for $Y$ and $\ov Y$ give a commutative diagram
of abelian groups with exact rows 
\begin{equation}
\begin{array}{ccccccccc}
0&\to &(\NS(\ov Y)/\ell^m)^\Ga &\to& \H^2(\ov Y,\mu_{\ell^m})^\Ga&\to& 
\Br(\ov Y)[\ell^m]^\Ga&&\\
&&\uparrow&&\uparrow&&\uparrow&&\\
0&\to & \Pic(Y)/\ell^m &\to& \H^2(Y,\mu_{\ell^m})&\to& 
\Br(Y)[\ell^m]&\to&0
\label{d1}
\end{array}
\end{equation}
If $(\NS(\ov Y)/\ell^m)^\Ga \to\H^2(\ov Y,\mu_{\ell^m})^\Ga$
is an isomorphism, then 
$\H^2(\ov Y,\mu_{\ell^m})^\Ga\to \Br(\ov Y)[\ell^m]^\Ga$
is the zero map. In this case from the commutativity of the right hand square of (\ref{d1})
and the surjectivity of $\H^2(Y,\mu_{\ell^m})\to\Br(Y)[\ell^m]$
we see that $\Br(Y)[\ell^m]\to \Br(\ov Y)[\ell^m]^\Ga$
is the zero map.
This shows that $\Br(Y)[\ell^m]$ is contained in $\Br_1(Y)$ for any $m$,
hence $\Br(Y)\{\ell\}\subset\Br_1(Y)$. In the particular case $Y=A$ 
we get $\Br(A)\{\ell\}\subset\Br_1(A)$.

The variety $Y$ is obtained by twisting $A$ by a cocycle
with coefficients in $A[2]$ acting on $A$ by translations.
The argument used in the proof of Proposition \ref{Br} shows that
the divisible group $A(\bar k)$, which contains $A[2]$, acts trivially on
the finite group $\H^2(\ov A,\mu_{\ell^m})$. 
Since $\NS(\ov Y)$ is canonically isomorphic
to $\NS(\ov A)$ as a $\Ga$-module, we have an isomorphism
of $\Ga$-modules $\NS(\ov Y)/\ell^m\cong \NS(\ov A)/\ell^m$
compatible with the cycle class map to 
$\H^2(\ov Y,\mu_{\ell^m})\cong\H^2(\ov A,\mu_{\ell^m})$. Thus
the injective map $(\NS(\ov Y)/\ell^m)^\Ga \to\H^2(\ov Y,\mu_{\ell^m})^\Ga$
is the same as $(\NS(\ov A)/\ell^m)^\Ga \to\H^2(\ov A,\mu_{\ell^m})^\Ga$.
It remains to show that this last map is an isomorphism. 

We have canonical isomorphisms of $\Ga$-modules
$$\wedge^2_{\Z/\ell^m}(\oplus_{i=1}^n A_i[\ell^m])=
(\oplus_{i=1}^n \wedge^2_{\Z/\ell^m}A_i[\ell^m]) \oplus
(\oplus_{i<j} (A_i[\ell^m]\otimes_{\Z/\ell^m} A_j[\ell^m]))$$
and
$$\Hom(A_i[\ell^m]\otimes_{\Z/\ell^m} A_j[\ell^m],\mu_{\ell^m})=
\Hom(A_i[\ell^m], \Hom(A_j[\ell^m],\mu_{\ell^m})).$$
Since each $A_i$ is principally polarised, the $\Ga$-modules 
$\Hom(A_j[\ell^m],\mu_{\ell^m})=A_j^t[\ell^m]$ 
and $A_j[\ell^m]$ are isomorphic. Hence the $\Ga$-module
\begin{equation}
\H^2(\ov A,\mu_{\ell^m})=\wedge^2_{\Z/\ell^m} \H^1(\ov A,\Z/\ell^m)(1)=
\Hom(\wedge^2_{\Z/\ell^m} A[\ell^m],\mu_{\ell^m}) \label{eq3}
\end{equation}
is isomorphic to the $\Ga$-module
\begin{equation}
\bigoplus_{i=1}^n \Hom(\wedge^2_{\Z/\ell^m}(A_i[\ell^m]),\mu_{\ell^m}) 
\oplus \bigoplus_{i<j} \Hom(A_i[\ell^m], A_j[\ell^m]). \label{eq33}
\end{equation}
For $i\not=j$ the $\Ga$-modules $A_i[\ell]$ and $A_j[\ell]$ are simple
and non-isomorphic, hence $\Hom_\Ga(A_i[\ell], A_j[\ell])=0$. 
We claim that
$\Hom_{\Ga}(A_i[\ell^m],A_j[\ell^m])=0$ for any $m\geq 1$ when $i\not=j$.
For $m>1$ the exact sequence of $\Ga$-modules
$$0 \lra A_i[\ell] \lra A_i[\ell^m] \stackrel{[\ell]}\lra A_i[\ell^{m-1}] \lra 0,$$
gives rise to an exact sequence of $\Z/\ell^m$-modules
$$0 \lra \Hom_{\Ga}(A_i[\ell^{m-1}], A_j[\ell^m])
\lra \Hom_{\Ga}(A_i[\ell^m], A_j[\ell^m])
\lra \Hom_{\Ga}(A_i[\ell], A_j[\ell^m]).$$
It is clear that
$$\Hom_{\Ga}(A_i[\ell^{m-1}], A_j[\ell^m])=
\Hom_{\Ga}(A_i[\ell^{m-1}], A_j[\ell^{m-1}]),$$ and
$$\Hom_{\Ga}(A_i[\ell], A_j[\ell^m])=
\Hom_{\Ga}(A_i[\ell], A_j[\ell]).$$
We obtain an exact sequence
\begin{equation}
0\to\Hom_{\Ga}(A_i[\ell^{m-1}], A_j[\ell^{m-1}])
 \to \Hom_{\Ga}(A_i[\ell^m], A_j[\ell^m])
 \to \Hom_{\Ga}(A_i[\ell], A_j[\ell]).\label{eqq}
\end{equation}
The induction assumption now implies
$\Hom_{\Ga}(A_i[\ell^m], A_j[\ell^m])=0$ when $i \ne j$.

If $\dim(A_i)=1$, then 
$\Hom(\wedge^2_{\Z/\ell^m}(A_i[\ell^m]),\mu_{\ell^m})$ is the trivial
$\Ga$-module $\Z/\ell^m$.

Now assume $\dim(A_i)>1$. Since the $\Ga$-modules $\Hom(A_i[\ell^m],\mu_{\ell^m})$ 
and $A_i[\ell^m]$ are isomorphic, the $\Ga$-module 
$\Hom(\wedge^2_{\Z/\ell^m} A_i[\ell^m],\mu_{\ell^m})$
is a submodule of $\End(A_i[\ell^m])$. Since
the $\Ga$-module $A_i[\ell]$ is absolutely simple, we have 
$\End_\Ga(A_i[\ell])=\F_\ell\cdot\Id$. We claim that
$\End_\Ga(A_i[\ell^m])=\Z/\ell^m\cdot \Id$ for any $m\geq 1$.
We argue by induction in $m$ and assume that
$$\End_{\Ga}(A_i[\ell^{m-1}])=\Z/\ell^{m-1} \cdot \Id.$$
In particular, the order of $\End_{\Ga}(A_i[\ell^{m-1}])$ equals $\ell^{m-1}$. 
The exact sequence (\ref{eqq}) in the case $i=j$ implies that the order of
$\End_{\Ga}(A_i[\ell^m])$ divides $\ell^{m-1}\cdot \ell=\ell^{m}$.
However, $\End_{\Ga}(A_i[\ell^m])$
contains the subgroup $\Z/\ell^m \cdot \Id$ of order $\ell^m$. This implies that
$\End_{\Ga}(A_i[\ell^m])=\Z/\ell^m \cdot \Id$, which proves our claim.

From (\ref{eq3}) and  (\ref{eq33}) we now conclude that $\H^2(\ov A,\mu_{\ell^m})^\Ga\subset(\Z/\ell^m)^n$.

The principal polarisation of each $A_i$ defines a non-zero class
in $\NS(\ov A_i)^\Ga$. It is well known that the $\Ga$-module $\oplus_{i=1}^n\NS(\ov A_i)$
is a direct summand of $\NS(\ov A)$, see, e.g. \cite[Prop. 1.7]{SZ4}.
Hence $\NS(\ov A)$ contains the trivial $\Ga$-module
$\Z^n$ as a full sublattice. Thus $(\NS(\ov A)/\ell^m)^\Ga$ contains
a subgroup isomorphic to $(\Z/\ell^m)^n$. 
It follows that the map $(\NS(\ov A)/\ell^m)^\Ga \to\H^2(\ov A,\mu_{\ell^m})^\Ga$
is an isomorphism for any $m\geq 1$. 

We have proved that $\Br(Y)\{\ell\}\subset\Br_1(Y)$. An equivalent statement is
that the natural map $\Br(Y)\{\ell\}\to\Br(\ov Y)$ is zero.
In the characteristic zero case Remark 2 in Section \ref{geo} implies that the natural map
$\Br(X)\{\ell\}\to \Br(\ov X)$ is zero. Equivalently, $\Br(X)\{\ell\}\subset\Br_1(X)$.
$\Box$

\medskip

Under additional assumptions we can prove a bit more.

\bpr \label{uno}
Let $k$ be a field of characteristic $0$. Let $\ell$ be a prime.
Let $A_1,\ldots,A_n$ be principally polarised abelian varieties over $k$ 
satisfying the following conditions.

{\rm (a)} The fields $k(A_i[\ell])$, where $i=1,\ldots,n$, are linearly disjoint over $k$.

{\rm (b)} The $\Ga$-module $A_i[\ell]$ is absolutely simple for each $i=1,\ldots,n$. 

{\rm (c)} $\NS(\ov A_i)\cong\Z$ for each $i=1,\ldots,n$. 

{\rm (d)} For each $i=1,\ldots,n$ the group $\Gal(k(A_i[\ell])/k)$ contains a subgroup $H_i$
which has no normal subgroup of index $\ell$, such that the $H_i$-module
$A_i[\ell]$ is simple, and, moreover, the $H_i$-module
$A_i[\ell]$ is absolutely simple if $\dim(A_i)>1$.

\noindent Let $A=\prod_{i=1}^n A_i$. Then $\Br(\ov A)[\ell]^\Ga=0$.
When $\dim(A)\geq 2$, for the Kummer variety $X$ attached to a $2$-covering of
$A$ we have $\Br(\ov X)[\ell]^\Ga=0$.
\epr
{\em Proof.} We claim that $\NS(\ov A)\cong \oplus_{i=1}^n\NS(\ov A_i)$.
It is well known that this is equivalent to the condition 
$\Hom(\ov A_i,\ov A_j)=0$ for all $i\not=j$, see, e.g. \cite[Prop. 1.7]{SZ4}. 
In view of (a) and (b) this condition holds by \cite[Thm. 2.1]{Z2}.
Now (c) implies that $\NS(\ov A)$ is isomorphic to the trivial $\Ga$-module $\Z^n$.

The properties $\Br(\ov A)[\ell]^\Ga=0$ and $\Br(\ov X)[\ell]^\Ga=0$
can be proved over any extension $k'$ of $k$ contained in $\bar k$. 
Let $k'$ be the compositum of $k(A_i[\ell])^{H_i}$ for $i=1,\ldots,n$,
and let $H=\prod_{i=1}^n H_i$.
Then $\Gal(k'(A[\ell])/k')=H$ and the fields $k'(A_i[\ell])$ 
are linearly disjoint over $k'$. By assumption (d) each $A_i[\ell]$ is a simple
$\Gal(\bar k/k')$-module and is absolutely simple whenever $\dim(A_i)>1$. 
Thus the assumptions of Proposition \ref{duo} are satisfied 
for the abelian varieties $A_1,\ldots,A_n$ over $k'$.
In the rest of the proof we write $k$ for $k'$ and $\Ga$ for $\Gal(\bar k/k')$.

The Kummer sequence gives an exact sequence of $\Ga$-modules
\begin{equation}
0\lra \NS(\ov A)/\ell\lra \H^2(\ov A,\mu_\ell)\lra \Br(\ov A)[\ell]\lra 0.
\label{a2}
\end{equation}
In view of (\ref{eq3}), $\Ga$ acts on the terms of (\ref{a2}) via 
its quotient $H$.
In particular, $\Br(\ov A)[\ell]^\Ga=\Br(\ov A)[\ell]^H$.
We obtain an exact sequence of cohomology groups of $H$:
\begin{equation} \label{a6}
0\lra (\NS(\ov A)/\ell)^H\lra \H^2(\ov A,\mu_\ell)^H \lra \Br(\ov A)[\ell]^H \lra
\H^1(H,\NS(\ov A)/\ell).
\end{equation}
The proof of Proposition \ref{duo} shows that 
the second arrow in (\ref{a6}) is an isomorphism. 
Since $\NS(\ov A)/\ell$ is the trivial $H$-module $(\F_\ell)^n$, we have
$$\H^1(H,\NS(\ov A)/\ell)=\Hom(H,(\F_\ell)^n)=0,$$
because by assumption $H$ has no normal subgroup of index $\ell$.
We conclude that
$\Br(\ov A)[\ell]^\Ga=\Br(\ov A)[\ell]^H=0$.
The second claim follows from Proposition \ref{Br}. $\Box$ 

\medskip

Note that the condition that $H$ has no normal
subgroup of index $\ell$ cannot be removed. See the remark on
\cite[p. 20]{SZ2} for an example of an abelian surface $A$ 
for which a Galois-invariant element in
$\Br(\ov A)[2]$ does not come from a Galois-invariant
element of $\H^2(\ov A,\mu_2)$. (In this example $H=\GL(2,\F_2)=\bS_3$,
the symmetric group on three letters.)

\medskip

Here is one of the main results of this paper.

\bthe \label{t1}
Let $k$ be a field of characteristic zero.
Let $A_1,\ldots,A_n$ be principally polarised abelian varieties over $k$ such that 
for $i=1,\ldots,n$ we have $\NS(\ov A_i)\cong\Z$, the $\Ga$-modules
$A_i[2]$ are absolutely simple, the fields $k(A_i[2])$
are linearly disjoint over $k$, and $\H^1(G_i,A_i[2])=0$, where
$G_i=\Gal(k(A_i[2])/k)$. Let $A=\prod_{i=1}^n A_i$. If $g=\dim(A)\geq 2$, then
for the Kummer variety $X$ attached to any $2$-covering of $A$
we have the following isomorphisms of abelian groups:

{\rm (i)} $\Pic(\ov X)\cong\Z^{2^{2g}+n}$;

{\rm (ii)} $\Br(X)\{2\}=\Br_1(X)\{2\}$;

{\rm (iii)} $\Br_1(X)=\Br_0(X)$.
\ethe 
{\em Proof.} Let $Y$ be the 2-covering of $A$ to which $X$ is attached.
It is clear that $Y=\prod_{i=1}^nY_i$, where $Y_i$ is a 2-covering of $A_i$
for $i=1,\ldots,n$. Using the principal polarisation of $A_i$ we identify $A_i$
with its dual abelian variety $A_i^t$. By \cite[Thm. 2.1]{Z2} 
we have $\Hom(\ov A_i,\ov A_j)=0$ for any $i\not=j$. A well known
consequence of this (see e.g., \cite[Prop. 1.7]{SZ4}) gives canonical
isomorphisms of $\Ga$-modules 
\begin{equation}
\Pic(\ov Y)\cong \bigoplus_{i=1}^n\Pic(\ov Y_i),\quad
\NS(\ov Y)\cong \bigoplus_{i=1}^n\NS(\ov Y_i)\cong\Z^n, \label{eq1}
\end{equation}
where $\Ga$ acts trivially on $\Z^n$. Now (i) follows from
the exact sequence (\ref{a1}).

Part (ii) follows from Proposition \ref{duo}, so it remains to establish part (iii).

For each $i=1,\ldots,n$ we have $Y_i=(A_i\times_k T_i)/A_i[2]$, where $T_i$ is
a torsor for $A_i[2]$, so that $Y=(A\times_k T)/A[2]$ with $T=\prod_{i=1}^n T_i$. 
Write $K_i=k(A_i[2])$ for the field of definition of the $2$-torsion subgroup of $A_i$
so that $G_i=\Gal(K_i/k)$.
Let $k(T_i)$ be the smallest subfield of $\bar k$
over which all $\bar k$-points of $T_i$ are defined.
Then $\Ga$ acts on $T_i(\bar k)\cong A_i[2]$ through $\Gal(k(T_i)/k)$.

If $T_i$ is a trivial torsor, then $T_i\cong A_i[2]$ and $\Gal(k(T_i)/k)=G_i$.
If $T_i$ is a non-trivial torsor, in our assumptions 
\cite[Prop. 3.6]{HS} gives us that $\Gal(k(T_i)/k)=
A_i[2]\rtimes G_i$, where $A_i[2]$ acts on itself by translations and 
$G_i$ acts on $A_i[2]$ by linear transformations.

The following lemma is a version of \cite[Prop. 3.12]{HS}.

\ble \label{9oct}
In the assumptions of Theorem \ref{t1}
the Galois extensions $k(T_1),\ldots,k(T_n)$ of $k$ are linearly disjoint over~$k$.
\ele
{\em Proof.} Let $m\leq n$ be the cardinality of $I\subseteq\{1,\ldots,n\}$
such that $T_i$ is a non-trivial torsor if and only if $i\in I$.
We proceed by double induction in $n\geq 1$ and $m\geq 0$. The statement
holds when $n=1$ (trivially) or when $m=0$ (by assumption).
Suppose that $n\geq 2$ and $m\geq 1$ and
the statement is proved for $(n,m-1)$ and for $(n-1,m-1)$.
Without loss of generality we can assume that $T_n$ is non-trivial.
By inductive assumption for $(n-1,m-1)$
the fields $k(T_i)$, $i=1,\ldots,n-1$, are linearly disjoint over $k$.
Let $L$ be the compositum of these fields, and let $E=L\cap k(T_n)$.
Each field $k(T_i)$ is Galois over $k$.
To check our statement it is enough to show that $E=k$.
By \cite[Cor. 3.9]{HS} the fact that $T_n$ is non-trivial implies that $E\subset K_n$
or $K_n\subset E$. By inductive assumption for $(n,m-1)$ we have $L\cap K_n=k$.
Thus $E\subset K_n$ implies $E=k$. On the other hand, $K_n\subset E$ implies $K_n=k$,
which is incompatible with our assumption that $A_n[2]$ is a simple $\Ga$-module. $\Box$

\medskip

Without loss of generality we can assume that $T_i$ is non-trivial for
$i=1,\ldots,m$ and $T_i$ is trivial for $i=m+1,\ldots,n$. Lemma \ref{9oct}
implies that the image of the action of $\Ga$ on $T(\bar k)=\prod_{i=1}^n T_i(\bar k)$
is the direct product
$$P=\prod_{i=1}^m (A_i[2]\rtimes G_i) \times \prod_{i=m+1}^n G_i.$$
Write $G=\prod_{i=1}^n G_i$. 
If we define $B=\prod_{i=1}^m A_i$, then $P=B[2]\rtimes G$.

To prove the desired property $\Br_1(X)=\Br_0(X)$ it is enough
to prove that $\H^1(k,\Pic(\ov X))=0$. 
The abelian groups in the exact sequence (\ref{a1}) are torsion-free, hence
$$\Pic(\ov X)\subset\Pic(\ov X)\otimes\Q\cong \Q[T]\oplus (\NS(\ov Y)\otimes\Q)
\cong \Q[T]\oplus\Q^n,$$
where $\Q[T]$ is the vector space with basis $T(\bar k)$ and
a natural action of $\Ga$. It follows that the image of the action of $\Ga$
on $\Pic(\ov X)$ is $P$. Thus it is enough to prove that
\begin{equation}
\H^1(P, \Pic(\ov X))=0. \label{c1}
\end{equation}

As an abelian group, $\Pi_1$ is generated by $\Z[T]$
and one half of the sum of the canonical generators of $\Z[T]$.
This gives an exact sequence of $A[2]\rtimes G$-modules
\begin{equation}
0\lra \Z[T]\lra\Pi_1\lra\Z/2\lra 0.\label{n2}
\end{equation}
By Shapiro's lemma $\H^1(P,\Z[T])=0$, because $\Z[T]$ is a permutation $P$-module. 
The cohomology exact sequences of (\ref{n2})
considered with respect to the action of $P$ and $G$
give rise to the following commutative diagram with exact upper row,
where the vertical arrows are given by restriction to the subgroup
$G\subset P$:
\begin{equation}\xymatrix{
0\ar[r]&\H^1(P,\Pi_1)\ar[r]\ar[d]&\Hom(P,\Z/2)\ar[d]\\
&0=\H^1(G,\Pi_1)\ar[r]&\Hom(G,\Z/2)}\label{n1}
\end{equation}
Here $\H^1(G,\Pi_1)=0$, as $\Pi_1$ is a permutation $G$-module,
see Remark 1 in Section \ref{geo}. For $i=1,\ldots, m$ we see from
\cite[Lemma 3.2 (ii)]{HS} that any subgroup of $A_i[2]\rtimes G_i$ 
of index 2 has the form $A_i[2]\rtimes H$ for a subgroup $H\subset G_i$ of index 2. 
Hence the right vertical arrow in (\ref{n1}) is an isomorphism.
The commutativity of (\ref{n1}) now implies that $\H^1(P,\Pi_1)=0$. 
The exact sequence of the cohomology groups of $P$
defined by the middle column of (\ref{diag})
shows that to prove (\ref{c1}) it is enough
to prove that $\H^1(P, \Pic(\ov Y)^{\iota_Y})=0$.

In view of the decomposition (\ref{eq1}) for this we must show that
$\H^1(P, \Pic(\ov Y_i)^{\iota_Y})=0$ for each $i=1,\ldots,n$.
For $Y_i$ the exact sequence (\ref{b3}) takes the form
\begin{equation}
0\lra A_i^t[2]\lra \Pic(\ov Y_i)^{\iota_Y}\lra \NS(\ov A_i)\lra 0. \label{eq2}
\end{equation}
Since $\NS(\ov A_i)\cong\Z$ we have $\H^1(P,\NS(\ov A_i))=0$. 

We first consider the case when $T_i$ is a trivial torsor. 
By assumption $\H^1(G_i,A_i[2])=0$.
We have $A_i[2]^{G_i}=0$, because $A_i[2]$ is a simple $G_i$-module with a
non-trivial action of $G_i$.
The restriction-inflation sequence for the normal subgroup $G_i\subset P$ 
acting on $A_i[2]$ shows that $\H^1(P,A_i[2])=0$, hence
$\H^1(P,\Pic(\ov Y_i)^{\iota_Y})=0$.

Now suppose that the torsor $T_i$ is non-trivial.
The Galois group $\Ga$ acts on $A_i[2]$ via $G_i$, hence so does
$\Gal(k(T_i)/k)=A_i[2]\rtimes G_i$.
We have $\H^1(A_i[2]\rtimes G_i,A_i[2])=\F_2$, see \cite[Prop. 3.6]{HS}.
This group is naturally a subgroup of $\H^1(k,A_i[2])$ and contains
the class $[T_i]$, because this class goes to zero under the restriction map
$\H^1(k,A_i[2])\to \H^1(k(T_i),A_i[2])$. Thus $[T_i]$ is the unique non-zero element of 
$\H^1(A_i[2]\rtimes G_i,A_i[2])$. 

Using the fact that $A_i[2]^{G_i}=0$, the Hochschild--Serre spectral sequence
for the normal subgroup $A_i[2]\rtimes G_i\subset P$ gives $\H^1(P,A_i[2])=\F_2$.
The same argument as above shows that $[T_i]$ is the unique non-zero element of
this group.

The principal polarisation $\lambda\in\NS(\ov A_i)^\Ga=\NS(\ov A_i)$ gives rise to
an isomorphism $\varphi_\lambda:A_i\tilde\lra A_i^t$ which
induces an isomorphism of $\Ga$-modules 
$\varphi_{\lambda*}:A_i[2]\tilde\lra A_i^t[2]$.
Since the $\Ga$-modules $\NS(\ov Y_i)$ and $\NS(\ov A_i)$
are canonically isomorphic, we can think of $\lambda$
as a generator of the trivial $\Ga$-module $\NS(\ov Y_i)\cong\Z$. 

Consider the exact sequence (\ref{eq2}) as a sequence of $P$-modules.
We claim that the differential $\NS(\ov Y_i)\to \H^1(P,A_i^t[2])$
sends the principal polarisation $\lambda$ to $\varphi_{\lambda*}[T_i]$, so this differential is surjective.
This implies that the first map in the exact sequence
$$
\H^1(P,A_i^t[2])\lra \H^1(P, \Pic(\ov Y_i)^\iota)
\lra \H^1(P,\Z)=0 $$
is zero, hence $\H^1(P, \Pic(\ov Y_i)^\iota)=0$.

To finish the proof of the theorem it remains to justify our claim. 
In the particular case of a trivial 2-covering 
the exact sequence of $\Ga$-modules (\ref{b3}) takes the form
\begin{equation}
0\lra A_i^t[2]\lra\Pic(\ov A_i)^{[-1]^*}\lra\NS(\ov A_i)\lra 0.\label{n4}
\end{equation}
Following \cite{PR} we shall write $c_\lambda$ for the image of 
$\lambda$ under the differential $\NS(\ov A_i)^{\Ga}
\to \H^1(k,A_i^t[2])$ attached to (\ref{n4}). By \cite[Lemma 3.6 (a)]{PR} we know that
$c_\lambda$ belongs to the kernel of the restriction map 
$\H^1(k,A_i^t[2])\to \H^1(K_i,A^t[2])$. We have $G_i=\Gal(K_i/k)$, so the restriction-inflation
sequence shows that $c_\lambda$ belongs to the subgroup
$\H^1(G_i,A_i^t[2])$ of $\H^1(k,A_i^t[2])$. However, our assumptions imply that
this group is zero, hence $c_\lambda=0$. Thus $\lambda$ is the image of some
$\Ga$-invariant
$\L\in\Pic(\ov A_i)^{[-1]^*}$, hence (\ref{n4}) is a {\em split} exact sequence
of $\Ga$-modules.

The exact sequence (\ref{eq2}) is obtained by twisting the exact sequence (\ref{n4})
by a 1-cocycle $\tau:\Ga\to A_i[2]$ representing $[T_i]\in\H^1(k,A_i[2])$.
By the definition of $\varphi_{\lambda}$ the translation by $x\in A_i(\bar k)$ acts
on $\Pic(\ov A_i)$ by sending $y\in\Pic(\ov A_i)$ to $y+\varphi_{\lambda}(x)$.
Since $Y_i$ is the twist of $A_i$ by $\tau$ with respect to the action of $A_i[2]$
by translations, we see that $g\in\Ga$ acts on $\L$, understood as an element
of $\Pic(\ov Y_i)$, by sending it to $\L+\varphi_{\lambda*}(\tau(g))$. By a standard
explicit description of the differential $\NS(\ov Y_i)\to \H^1(k,A_i^t[2])$
we see that $\lambda$ goes to $\varphi_{\lambda*}[T_i]$, as claimed. 
$\Box$

\medskip

\noindent{\bf Remark 4}
If $k$ is finitely generated over $\Q$, then $\Br(X)/\Br_0(X)$ is finite
by Corollary \ref{cor2}. 
The 2-primary subgroup of $\Br(X)/\Br_0(X)$ is the image of
the 2-primary subgroup of $\Br(X)$, and hence Theorem \ref{t1} implies that
the order of $\Br(X)/\Br_0(X)$ is odd.

\bco \label{HaSk}
Let $k$ be a number field and let $X$ be a Kummer variety
satisfying the assumptions of Theorem \ref{t1}. If $X(\A_k)\not=\emptyset$,
then $X(\A_k)^\Br\not=\emptyset$.
\eco
{\em Proof.} In view of Remark 4 this is a formal consequence of Theorems \ref{g1} and \ref{t1}. $\Box$

\medskip

This corollary explains the absence of the Brauer--Manin obstruction from the 
statement of the (conditional) Hasse principle for Kummer varieties recently established
in \cite[Thm. 2.2]{HS}, once we impose the additional condition $\NS(\ov A_i)\cong\Z$ for
$i=1,\ldots,n$. In the next section we give examples related to hyperelliptic curves
where this condition holds.

\section{Kummer varieties attached to products of Jacobians of hyperelliptic curves}

In this section we consider the case when each factor of $A=\prod_{i=1}^nA_i$ is the Jacobian of
a hyperelliptic curve given by a polynomial of odd degree $\geq 3$ with a large Galois group.
It will be convenient to include elliptic curves as a particular case of hyperelliptic curves,
so we shall adopt this terminology here without further mention.

We write $\bS_n$ for the symmetric group on $n$ letters,
and $\bA_n\subset\bS_n$ for the alternating group on $n$ letters.

\bthe \label{t11}
Let $k$ be a field of characteristic zero.
Let $A$ be the product of Jacobians
of the hyperelliptic curves $y^2=f_i(x)$, where $f_i(x)\in k[x]$ is a
separable polynomial of odd degree $d_i\geq 5$ 
with Galois group $\bS_{d_i}$ or $\bA_{d_i}$, 
or of degree $3$ with Galois group $\bS_3$, for $i=1,\ldots,n$.
Assume that $g=\sum_{i=1}^n(d_i-1)/2\geq 2$ and 
the splitting fields of the polynomials $f_i(x)$, $i=1,\ldots,n$,
are linearly disjoint over $k$.
Then the conclusions of Theorem \ref{t1} hold for the Kummer variety $X$
attached to any $2$-covering of $A$. Moreover, $\Br(\ov X)[2]^\Ga=0$.
\ethe
{\em Proof.} For $i=1,\ldots,n$ let $C_i$ be the smooth and projective curve given by
the equation $y^2=f_i(x)$ and let $A_i$ be the Jacobian of $C_i$.
Let $A=\prod_{i=1}^nA_i$, and let $Y$ be a 2-covering of $A$ such that $X$ is the 
Kummer variety attached to $Y$.

Since $A_i$ is canonically principally polarised,
we have an isomorphism $A_i\tilde\lra A_i^t$. It is well known that
$\NS(\ov A_i)$ is isomorphic to the subgroup of self-dual endomorphisms
in $\End(\ov A_i)=\Hom(\ov A_i,\ov A_i^t)$. If $\deg(d_i)\geq 5$,
by \cite[Thm. 2.1]{Z} we have $\End(\ov A_i)\cong\Z$, hence $\NS(\ov A_i)\cong\Z$. 
If $\deg(d_i)=3$, then we obviously have $\NS(\ov A_i)\cong\Z$.

Let $W_i\subset C_i$ be the subscheme given by $f_i(x)=0$.
The double covering $C_i\to\P^1_k$ is ramified precisely
at the $\bar k$-points of $W_i\cup\{\infty\}$. 
It is well known that
the $\Ga$-module $A_i[2]$ is isomorphic to the zero sum subspace of
the $\F_2$-vector space with basis $W_i(\bar k)$, with the action of $\Ga$
defined by the natural action of $\Ga$ on $W(\bar k)$.
In particular, the splitting field of $f_i(x)$ is $k(A_i[2])$,
and the Galois group of $f_i(x)$ is $G_i=\Gal(k(A_i[2])/k)$.
This implies that the fields $k(A_i[2])$ are linearly disjoint over $k$.
Since $d_i$ is odd, the permutation $G_i$-module $\F_2^{d_i}$ whose 
canonical generators
are given by the $\bar k$-points of $W$, is the direct sum $\F_2\oplus A_i[2]$.
We note that for $d_i\geq 5$ the standard representation of 
$\bA_{d_i}\subset {\rm Sp}(d_i-1,\F_2)$ in $\F_2^{d_i-1}$ 
is absolutely irreducible \cite{M} (see also \cite[Lemma 5.2]{Z}).
The same is true if we replace $\bA_{d_i}$ by $\bS_{d_i}$.
If $d_i=3$, then the standard 2-dimensional representation of $\bS_3$ is
absolutely irreducible \cite{M}.
This implies that in all our cases $\End_{G_i}(A_i[2])=\F_2$,
cf. \cite[Thm. 5.3]{Z}.
We conclude that the $\Ga$-module $A_i[2]$ is absolutely simple, for $i=1,\ldots,n$.

By Shapiro's lemma we have
$\H^1(\bA_{d_i},\F_2[\bA_{d_i}/\bA_{d_i-1}])=\H^1(\bA_{d_i-1},\F_2)=0$ 
for $d_i\geq 5$,
because $\bA_{d_i-1}$ has no subgroup of index 2 (it is generated by elements of order 3).
Hence $\H^1(\bA_{d_i},A_i[2])=0$. Similarly, 
$\H^1(\bS_{d_i},\F_2[\bS_{d_i}/\bS_{d_i-1}])=\H^1(\bS_{d_i-1},\F_2)=\F_2$
for $d_i\geq 3$,
because $\bA_{d_i-1}$ is the unique subgroup of $\bS_{d_i-1}$ of index 2.
This implies $\H^1(\bS_{d_i},A_i[2])=0$ (cf. \cite[Lemma 2.1]{HS}). Thus $\H^1(G_i,A_i[2])=0$
for all $i=1,\ldots,n$.

We have checked that all the assumptions of Theorem \ref{t1} are satisfied.
In particular, conditions (a), (b), (c) of Proposition \ref{uno} 
are satisfied. Condition (d) is also satisfied if we take $H_i=\bA_{d_i}$
for $i=1,\ldots,n$. Indeed, each $A_i[2]$ is a simple $\bA_{d_i}$-module for all odd
$d_i\geq 3$ and is absolutely simple if $d_i\geq 5$. Finally,
$\bA_{d_i}$ has no subgroup of index 2 as it is generated by the elements of order 3.
An application of Proposition \ref{uno} gives that $\Br(\ov X)[2]^\Ga=0$.
$\Box$

\bco \label{g2}
Let $k$ be a number field.
Let $A$ be the product of Jacobians
of the hyperelliptic curves $y^2=f_i(x)$, where $f_i(x)\in k[x]$ is a
separable polynomial of odd degree $d_i\geq 5$ 
with Galois group $\bS_{d_i}$ or $\bA_{d_i}$, 
or of degree $3$ with Galois group $\bS_3$, for $i=1,\ldots,n$.
Assume that $g=\sum_{i=1}^n(d_i-1)/2\geq 2$ and 
the splitting fields of the polynomials $f_i(x)$, $i=1,\ldots,n$,
are linearly disjoint over $k$.
If the Kummer variety $X$
attached to a $2$-covering of $A$ is everywhere locally soluble, then $X(\A_k)^\Br\not=\emptyset$.
\eco
{\em Proof.} In view of Remark 4 at the end of Section \ref{kuku}
this is a formal consequence of Theorems \ref{g1} and \ref{t11}. $\Box$

\bigskip

\noindent{\bf Example 1} L. Dieulefait shows in \cite[Thm. 5.8]{D}
that for $k=\Q$ and $f(x)=x^5-x+1$ the image of the Galois group
$\Ga=\Gal(\ov\Q/\Q)$ in $\Aut(A[\ell])$, where $A$ is 
the Jacobian of the hyperelliptic curve $y^2=f(x)$, is $\mathrm{GSp}(4,\F_\ell)$
for each prime $\ell\geq 3$. (In \cite{D} this result was conditional
on the Serre conjectures \cite{Serre}, which have been later proved by C. Khare
and J.-P. Wintenberger \cite{KW}.)
A verification with {\tt magma} gives that
the Galois group of $x^5-x+1$ is $\bS_5$, so Theorem \ref{t11}
can be applied. Thus for the Kummer surface $X$
attached to a 2-covering of $A$ we have $\Br(\ov X)[2]^\Ga=0$ and $\Br_1(X)=\Br_0(X)$.
On the other hand, Proposition \ref{uno} can be applied for 
each prime $\ell\geq 3$ with $H=\mathrm{Sp}(4,\F_\ell)$.
Indeed, for $\ell\geq 3$ the group $\mathrm{PSp}(4,\F_\ell)$ is simple non-abelian
\cite[Thm. 5.2, p. 177]{GA} of order $\ell^4(\ell^4-1)(\ell^2-1)/2>\ell$, so $H$
contains no normal subgroups of index $\ell$.
The tautological representation of $\mathrm{Sp}(4,\F_\ell)$
is well known to be absolutely irreducible. We obtain
that the Kummer surface $X$ attached to a 2-covering of $A$ 
is a K3 surface of geometric 
Picard rank 17 such that $\Br(\ov X)^\Ga=0$. Hence $\Br(X)=\Br_0(X)$. 

\medskip

\noindent{\bf Example 2} R. Jones and J. Rouse consider the Jacobian $A$ of the curve
of genus 2 given by $y^2=f(x)$, where $f(x)=4 x^6-8 x^5+4 x^4+4 x^2-8x+5$
is a polynomial over $\Q$ with Galois group $\bS_6$ and discriminant quadratic extension 
$\Q(\sqrt{-3\cdot 13\cdot 31})$, see \cite[Example 6.4, pp. 787--788]{JR}.  
They show that the image of $\Ga=\Gal(\ov\Q/\Q)$ in
$\Aut(A[\ell])$ is $\mathrm{GSp}(4,\F_\ell)$ for all primes $\ell$.
For odd $\ell$ the only non-trivial isomorphic quotients of 
$\mathrm{GSp}(4,\F_2)\cong \bS_6$ and 
$\mathrm{GSp}(4,\F_\ell)$ are cyclic groups of order 2, namely
$\bS_6/\bA_6$ and $\mathrm{GSp}(4,\F_\ell)/\F_{\ell}^{*2}\cdot\mathrm{Sp}(4,\F_\ell)$, 
respectively. By Goursat's lemma a subgroup of 
$\bS_6\times \mathrm{GSp}(4,\F_\ell)$ that maps surjectively onto each factor 
is either the whole product or the inverse image of the graph of the unique isomorphism 
$$\bS_6/\bA_6\ \tilde\lra\ \mathrm{GSp}(4,\F_\ell)/\F_{\ell}^{*2}\cdot\mathrm{Sp}(4,\F_\ell).$$
Hence such a subgroup contains $\bA_6 \times \mathrm{Sp}(4,\F_\ell)$. 
Let $\alpha$ be a root of $f(x)$, and let
$k=\Q(\alpha)$ or $k=\Q(\alpha, \sqrt{-3\cdot 13\cdot 31})$.
Then the Galois group of $f(x)$ over $k$ is $\bS_5$ or $\bA_5$, respectively,
whereas $\Gal(k(A[\ell])/k)$ contains $\mathrm{Sp}(4,\F_\ell)$ for all $\ell\geq 3$.
Now the same arguments as in Example 1 show that for the Kummer surface $X$ over $k$
attached to a 2-covering of $A$ 
we have $\Br(\ov X)^\Ga=0$ and $\Br_1(X)=\Br_0(X)$, hence $\Br(X)=\Br_0(X)$.

\medskip

\noindent{\bf Example 3} D. Zywina \cite[Thm. 1.1]{Zy} gives an example
of a smooth plane quartic curve over $\Q$ such that the image of 
$\Ga=\Gal(\ov\Q/\Q)$ on the torsion points of its Jacobian $A$ is
the full group $\mathrm{GSp}(6,\hat\Z)$. 
We have $\End(\ov A)\cong\Z$, as follows from \cite[Thm. 3, p. 577]{Z10},
where one takes $X=A$, $\tilde{G}_2=\mathrm{GSp}(6,\F_2)$ 
and $G=\mathrm{Sp}(6,\F_2)$. This implies $\NS(\ov A)\cong\Z$.
Let $k\subset \Q(A[2])$ be such that $\Gal(\Q(A[2])/k)$ is 
$\bS_7$ or $\bA_7$ embedded into $\mathrm{Sp}(6,\F_2)$ in the usual way.
We can adapt the proof of Theorem \ref{t11} to this case and use
Proposition \ref{uno} in the same way as in Example 1. This
shows that the Kummer threefold $X$ over $k$ attached to a 2-covering of $A$ 
has $\Br(\ov X)^\Ga=0$ and $\Br_1(X)=\Br_0(X)$, hence $\Br(X)=\Br_0(X)$.

\bigskip

\noindent Department of Mathematics, South Kensington Campus,
Imperial College London, SW7 2BZ England, U.K. -- and --
Institute for the Information Transmission Problems,
Russian Academy of Sciences, 19 Bolshoi Karetnyi, Moscow, 127994
Russia

\medskip

\noindent a.skorobogatov@imperial.ac.uk

\bigskip

\noindent Department of Mathematics, Pennsylvania State University, University Park, Pennsylvania 16802, USA 

\medskip

\noindent zarhin@math.psu.edu

\end{document}